\documentclass{scrartcl}

\usepackage{graphicx} 
\usepackage{amsmath, amssymb, amsfonts, amsthm, bm}
\usepackage{hyperref}
\usepackage{graphicx}
\usepackage{sidecap}
\usepackage[dvipsnames]{xcolor}
\colorlet{RefColor}{green!50!black}
\colorlet{LinkColor}{red!50!black}

\usepackage{algorithm}
\usepackage{algpseudocode}
\usepackage{mathtools}
\usepackage{authblk}

\usepackage[most]{tcolorbox}
\hypersetup{
  colorlinks = true,
  linkbordercolor = {white},
  linkcolor = LinkColor,
  anchorcolor=LinkColor,
  citecolor= RefColor,
  filecolor=cyan,
  menucolor=red,
  runcolor=cyan,
  urlcolor = LinkColor,
}
\usepackage{Notation}
\usepackage[font=small,labelfont=bf, format=plain]{caption}

\definecolor{lightgray}{HTML}{f3f3f3}

\newtheorem{remark}{Remark}
\newtheorem{proposition}{Proposition}

\title{Nonlinear embeddings for conserving Hamiltonians and other quantities with Neural Galerkin schemes\thanks{
    The authors Schwerdtner, Berman, Peherstorfer were partially supported by National Science Foundation under Grant No.~2046521 and the Office of Naval Research under award N00014-22-1-2728.}}

\author[1]{Paul Schwerdtner\thanks{corresponding author: paul.schwerdtner@nyu.edu}}
\author[2]{Philipp Schulze}
\author[1]{Jules Berman}
\author[1]{Benjamin Peherstorfer}
\affil[1]{Courant Institute of Mathematical Sciences, New York University}
\affil[2]{Technische Universit\"at Berlin}
\date{October 2023}

\usepackage[%
  defernumbers=true,
  maxbibnames=999,
  sortcites=true,
  giveninits=true,                 
  style=trad-plain
]{biblatex}

\newbibmacro{string+doiurl}[1]{%
  \iffieldundef{doi}{%
    \iffieldundef{url}{%
      #1%
    }{%
      \href{\thefield{url}}{#1}%
    }%
  }{%
    \href{https://doi.org/\thefield{doi}}{#1}%
  }%
}

\DeclareFieldFormat[book,thesis,incollection,inbook]{title}{\usebibmacro{string+doiurl}{\mkbibemph{#1}}}
\DeclareFieldFormat[article,inproceedings]{title}{\usebibmacro{string+doiurl}{#1}} 

\usepackage{tikz, pgfplots}
\pgfplotsset{compat=1.18}
\usetikzlibrary{arrows.meta}
\usetikzlibrary{backgrounds}
\usetikzlibrary{intersections,shapes.arrows}

\usepgfplotslibrary{patchplots}
\usepgfplotslibrary{fillbetween}
\pgfplotsset{%
  layers/standard/.define layer set={%
      background,axis background,axis grid,axis ticks,axis lines,axis tick labels,pre main,main,axis descriptions,axis foreground%
    }{
      grid style={/pgfplots/on layer=axis grid},%
      tick style={/pgfplots/on layer=axis ticks},%
      axis line style={/pgfplots/on layer=axis lines},%
      label style={/pgfplots/on layer=axis descriptions},%
      legend style={/pgfplots/on layer=axis descriptions},%
      title style={/pgfplots/on layer=axis descriptions},%
      colorbar style={/pgfplots/on layer=axis descriptions},%
      ticklabel style={/pgfplots/on layer=axis tick labels},%
      axis background@ style={/pgfplots/on layer=axis background},%
      3d box foreground style={/pgfplots/on layer=axis foreground},%
    },
}

\allowdisplaybreaks
\begin{document}

\maketitle

\begin{abstract}
	This work focuses on the conservation of quantities such as Hamiltonians, mass, and momentum when solution fields of partial differential equations are approximated with nonlinear parametrizations such as deep networks. The proposed approach builds on Neural Galerkin schemes that are based on the Dirac--Frenkel variational principle to train nonlinear parametrizations sequentially in time. We first show that only adding constraints that aim to conserve quantities in continuous time can be insufficient because the nonlinear dependence on the parameters implies that even quantities that are linear in the solution fields become nonlinear in the parameters and thus are challenging to discretize in time. Instead, we propose Neural Galerkin schemes that compute at each time step an explicit embedding onto the manifold of nonlinearly parametrized solution fields to guarantee conservation of quantities. The embeddings can be combined with standard explicit and implicit time integration schemes. Numerical experiments demonstrate that the proposed approach conserves quantities up to machine precision.
\end{abstract}

\section{Introduction}
Preserving structure and conserving quantities such as Hamiltonians, mass, and momentum in numerical approximations of solution fields governed by time-dependent partial differential equations (PDEs) is important to guarantee physical consistency and help interpretation.
This work focuses on conservation of quantities when discretizing solution fields with nonlinear parametrizations such as deep neural networks.
The challenge is that the nonlinear dependence on the parameters means that quantities that are linear in the solution fields of the PDEs become nonlinear in the parameters, which leads to a loss of the linear vector-space structure that numerical methods traditionally build on.
We show that only adding constraints in continuous-time formulations to keep quantities constant over time in the nonlinearly parametrized solution fields is insufficient.
Instead, we propose to compute explicit embeddings onto the manifold of nonlinearly parametrized fields that conserve quantities.
Numerical experiments demonstrate that the proposed approach conserves quantities up to machine precision.

Preserving structure and conserving quantities is a mainstay in computational science and engineering. There is a range of works that aim to learn models from data while preserving structure, such as structure-preserving dynamic mode decomposition \cite{BaddooHMKB2023Physics-informed,MorandinNU2023Port-Hamiltonian}, operator inference \cite{SharmaWK2022Hamiltonian,P19ReProj,SawantKP2023Physics-informed,QKPW19LiftLearn,GruberT2023Canonical_arxiv}, methods that learn from frequency-domain data \cite{SchulzeU2016Data-driven,SchulzeUBG2018Data-driven,PontesDuffGB2022Data-driven,SchwerdtnerV2023SOBMOR,GoseaGW2023Structured}, and methods based on deep learning \cite{ChenZAB2020Symplectic,otness2021an,ZhangCJ2022Implicit}.
Another line of work aims to preserve structure in reduced models \cite{PengM2016Symplectic,doi:10.1137/15M1055085,BuchfinkBH2019Symplectic,GongWW2017Structure-preserving,AfkhamH2017Structure,FarhatCA2015Structure-preserving,KalashnikovaBAB2014Stabilization}, specifically in the setting of computational fluid dynamics \cite{BaroneKST2009Stable,RowleyCM2004Model,Chan2020Entropy,BalajewiczTD2016Minimal,RosenbergerS2022Momentum-conserving,AfkhamRWH2020Conservative,BoonF2023Reduced,Liljegren-SailerM2022Snapshot-based,MohebujjamanRI2019Physically,Sanderse2020Non-linearly}.
Closest to our work is~\cite{CarlbergCS2018Conservative} that introduces a modified Galerkin system so that reduced solution fields conserve quantities. In all of these works, the parametrizations are linear.

We focus on nonlinear parametrizations, which, for example, can be trained with the Dirac--Frenkel variational principle \cite{Dirac1930Note,Frenkel1934Wave,Kramer1981,Lubich2008}. Under assumptions on the parametrization and the equation of interest, structure is preserved in the special case when the parametrization is so rich that the residual and thus the error vanishes, which is leveraged in, e.g., \cite{LasserL2020Computing}.
The work \cite{AndersonF2022Evolution} considers nonlinear parametrizations that are tuned towards the solution fields and proposes to add constraints that help conserve quantities in continuous time.
Another line of work proposes dynamic low-rank approximations that are structure-preserving such as \cite{EINKEMMER2023112060,EINKEMMER2021110353,Pagliantini2021Dynamical,HesthavenPR2022Rank-adaptive}.
In model reduction with nonlinear parametrizations \cite{P22AMS}, there is the work  \cite{BuchfinkGH2023Symplectic} that preserves symplecticity in a weak sense. The work \cite{SharmaMBGGK2023Symplectic_arxiv} is a major step forward and considers quadratic manifolds to achieve an offline/online decoupling but currently still incurs high online costs because no empirical interpolation is considered.
The work closest to ours is \cite{LeeC2021Deep}, which considers deep autoencoders for parametrizing the latent states and adds constraints in the time-discrete formulation to conserve quantities; however, it requires solving a nonlinear and potentially non-convex optimization problem at each time step with the number of unknowns scaling with the latent state dimension and the number of conserved quantities. In fact, the costs of computing the optimization objective and gradients grow in the dimension of the ambient space rather than the dimension of the latent space. In contrast, our approach is applicable to explicit time integration schemes that lead to linear regression problems in each time step. We also have to solve a system of nonlinear equations at each time step but the number of unknowns grows with the number of conserved quantities only, which typically is orders of magnitude lower than the state dimension.
Additionally, we leverage the results of the work \cite{Schulze2023Structure-Preserving} that achieves structure preservation via weighted schemes and specific nonlinear parametrizations, instead of explicitly computing embeddings. We show that the same specific nonlinear parametrizations lead to structure preservation in our setting too.

We build on the Neural Galerkin scheme introduced in \cite{BrunaPV2022Neural}, which is based on the Dirac--Frenkel variational principle \cite{Dirac1930Note,Frenkel1934Wave,Kramer1981,Lubich2008} and applies to generic nonlinear parametrizations such as deep networks.
We first introduce Neural Galerkin schemes with constraints that conserve quantities in continuous time.
We also show that Hamiltonians can be conserved in continuous time with weighted Neural Galerkin schemes and specific nonlinear parametrizations that are separable, even without adding constraints.
While systems may conserve quantities in continuous time, the time discretization is delicate because the nonlinear parametrization means that the quantities depend nonlinearly on the parameters and thus no Runge-Kutta integration scheme can exist that conserves such quantities in general \cite[Theorem IV.3.3]{HairerWL2006Geometric}.
Instead, we use nonlinear embeddings to find approximations that are close to the Neural Galerkin solutions and at the same time conserve quantities in discrete time. Importantly, the nonlinear embeddings can be combined with explicit time integration schemes that can be more efficient in the context of nonlinear parametrizations than implicit schemes \cite{BrunaPV2022Neural}. The nonlinear embeddings follow \cite{HairerWL2006Geometric} and are numerically computed via an iterative scheme that applies generically to nonlinear parametrizations including deep networks but also to other nonlinear parametrizations used in, e.g., model reduction such as \cite{Schulze2023Structure-Preserving}. We stress that our approach conserves quantities but not necessarily structure such as symplecticity.
With Burgers', acoustic wave, and shallow water equations, we demonstrate that the proposed scheme conserves quantities such as mass, energy, and Hamiltonians up to machine precision in numerical experiments.

\section{Preliminaries}
\label{sec:ng_explanation}
We describe the setup of time-dependent PDEs with conserved quantities and Hamiltonians, discuss Neural Galerkin schemes based on the Dirac--Frenkel variational principle, and provide a problem formulation.

\subsection{Setup of time-dependent PDEs}
Consider the PDE
\begin{align}
	\label{eq:general_pde}
	\partial_t \vecu(t, \vecx) & =  \vecf(\vecx, \vecu(t, \cdot)) & (t, \vecx) \in \timedomain \times \spacedomain\,, \\
	\vecu(0, \vecx)            & =  \vecu_0(\vecx)                & \vecx \in \spacedomain\,,
\end{align}
with the solution field $\vecu : \timedomain \times \spacedomain \rightarrow \R^{\outputdim}$ on the spatial domain $\spacedomain \subseteq \R^{\spacedim}$. At each time $t \in \timedomain$, the solution field $\vecu(t, \cdot): \spacedomain \to \R^{\outputdim}$ is in a space $\mathcal{U}$ of functions that allow point-wise evaluations. In the following, the space $\mathcal{U}$ is a subspace of the space $\Lp{2}(\spacedomain)^m$ of square-integrable functions with $\outputdim$ outputs with respect to a fully supported measure $\nu$.
The right-hand side function $\vecf: \spacedomain \times \mathcal{U} \rightarrow \R^{\outputdim}$ can include partial derivatives of $\vecu$ in the spatial variable $\vecx$. The initial condition is $\vecu_0 \in \mathcal{U}$.
In the following, the boundary conditions for equation \eqref{eq:general_pde} are imposed by seeking solutions in the space $\mathcal{U}$ so that~\eqref{eq:general_pde} is well posed.

\subsection{Conserved quantities}
Consider now  quantities of the form
\begin{align}
	\label{eq:conserved_quantity_format}
	\cq_i: \mathcal{U} \rightarrow \R \quad \vecu(t, \cdot) \mapsto \int_\spacedomain \quantitykernel_i(\vecu(t, \cdot))(\vecx) \dd \nu(\vecx),
\end{align}
where $\quantitykernel_i: \mathcal{U} \rightarrow \Lp{1}(\spacedomain)$ is continuously differentiable for $i  = 1, \dots, \ncq$.
We categorize a quantity $\cq_i$ as linear, quadratic or nonlinear depending on whether $\quantitykernel_i$ is a linear, quadratic or nonlinear function in $\vecu$, respectively.
A quantity $\cq_i$ is called a conserved quantity of the solution field $\vecu$ of equation \eqref{eq:general_pde} if it remains constant in the sense $\cq_i(\vecu(t, \cdot)) = \cq_i(\vecu(0, \cdot))$ for all $t \in \timedomain$ \cite[Chapter IV]{HairerWL2006Geometric}.
If the integrals in the conserved quantities~\eqref{eq:conserved_quantity_format} cannot be computed analytically, we numerically estimate them via Monte Carlo from $\nsamplesmani$ samples as
\begin{equation}
	\label{eq:conserved_quantity_monte_carlo}
	\hat\cq_i: \mathcal{U} \rightarrow \R \quad \vecu(t, \cdot) \mapsto \frac{1}{\nsamplesmani}\sum_{s = 1}^{\nsamplesmani} \quantitykernel_i(\vecu(t, \cdot))(\vecxi_s)\,, \quad i = 1, \dots, \ncq\,,
\end{equation}
where $\vecxi_1, \dots, \vecxi_{\nsamplesmani} \in \spacedomain$. The following methodology extends in a straightforward way to other quadrature schemes.

\subsection{Hamiltonian systems}
An important example of a conserved quantity arises when the right-hand side $\vecf$ of \eqref{eq:general_pde} can be written in the Hamiltonian form,
\begin{equation}
	\label{eq:structuredRHS}
	\vecf(\cdot, \vecv) = \Joperator(\vecv)\frac{\delta\hamiltonian}{\delta\vecu}(\vecv)
\end{equation}
for all $t\in\timedomain$ and $\vecv\in\mathcal{U}$, with the Hamiltonian $\hamiltonian\colon \mathcal{U}\to\R$ and the interconnection operator $\Joperator\colon \mathcal{U}\to\Hom(\mathcal{U},\mathcal{Z})$, where $\Hom(\mathcal{U},\mathcal{Z})$ denotes the space of linear operators from $\mathcal{U}$ to $\mathcal{Z}$.
The space $\mathcal{Z}$ is a subspace of functions of $\Lp{2}(\spacedomain)^{\outputdim}$ that allow point-wise evaluations. Additionally, the variational derivative $\frac{\delta\hamiltonian}{\delta\vecu}$ attains values in $\mathcal{U}$.
In the following, we consider Hamiltonians $\hamiltonian$ that can be written as
\begin{equation}
	\label{eq:hamiltonian_as_integral}
	\hamiltonian(\vecv) = \int_{\spacedomain}h (\vecv)(\vecx) \dd \nu(\vecx)\,,
\end{equation}
with continuously differentiable $h\colon \mathcal{U} \rightarrow \Lp{1}(\spacedomain)$.
Analogous to the sampled conserved quantities \eqref{eq:conserved_quantity_format}, we also introduce the sampled Hamiltonian $\Hsampled$.
The interconnection operator $\Joperator$ is
pointwise skew-adjoint in the sense that
\begin{equation}
	\label{eq:skew-adjointness_J}
	\langle \vecq_1, \Joperator(\vecv) \vecq_2 \rangle_\nu = -\langle\Joperator(\vecv)\vecq_1, \vecq_2 \rangle_\nu
\end{equation}
holds for all $\vecv,\vecq_1,\vecq_2\in\mathcal{U}$, where $\ip{\cdot,\; \cdot}_\nu$ denotes the $\Lp{2}$ inner product corresponding to the measure $\nu$.
The pointwise skew-adjoint property implies in particular
\begin{equation}
	\label{eq:vanishing_quadratic_form_associated_with_J}
	\langle \vecq, \Joperator(\vecv) \vecq \rangle_\nu = \frac12(\langle \vecq, \Joperator(\vecv) \vecq \rangle_\nu-\langle \Joperator(\vecv)\vecq, \vecq \rangle_\nu) = 0
\end{equation}
for all $\vecv,\vecq\in\mathcal{U}$.
The property \eqref{eq:vanishing_quadratic_form_associated_with_J} implies that the Hamiltonian is a conserved quantity, which follows from the computation
\begin{align*}
	\frac{\dd\hamiltonian(\vecu(t,\cdot))}{\dd t}(t)
	 & =
	\Big\langle \frac{\delta\hamiltonian}{\delta\vecu}(\vecu(t,\cdot)),\; \partial_t \vecu(t, \cdot) \Big\rangle_\nu
	\\
	 & =
	\Big\langle \frac{\delta\hamiltonian}{\delta\vecu}(\vecu(t,\cdot)),\; \Joperator(\vecu(t,\cdot))\frac{\delta\hamiltonian}{\delta\vecu}(\vecu(t,\cdot)) \Big\rangle_\nu  =0.
\end{align*}
One example of an equation that can be represented with the structure \eqref{eq:structuredRHS} is the inviscid Burgers' equation with periodic boundary conditions in one spatial dimension,
\begin{equation}
	\label{eq:burgers_equation_hamiltonian_structure}
	\Joperator(v)q = -\frac13\left(\partial_x(v q)+v\partial_x q\right),\quad \hamiltonian(u) = \frac12\lVert u\rVert_{\Lp{2}(\spacedomain)}^2\,,
\end{equation}
which we will discuss in Section~\ref{sec:burgers_equation}.
Another example is the linear wave equation with periodic boundary conditions,
\begin{equation}
	\label{eq:wave_equation_hamiltonian_structure}
	\Joperator = -
	\begin{bmatrix}
		0          & \partial_x \\
		\partial_x & 0
	\end{bmatrix}
	,\quad \hamiltonian(\rho,v) = \frac{1}{2}\int_{\spacedomain} \frac{c^2}{\rhoref}\rho(\vecx)^2+\rhoref v(\vecx)^2 \; \dd \vecx\,,
\end{equation}
which we will discuss in Section~\ref{sec:wave_equation}.

\subsection{Neural Galerkin schemes based on the Dirac--Frenkel variational principle}\label{sec:NG}
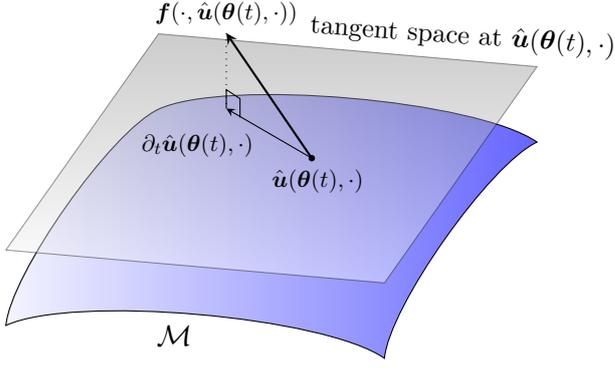
\begin{SCfigure}
	\begin{tikzpicture}

\tikzset{x={(175:1cm)},y={(55:.7cm)},z={(90:1cm)}}

\shade[left color=blue!4,right color=blue!60,looseness=.5, draw=black]  (2.5,-2.5,-1) 
to[bend left] (2.5,2.5,-1)
to[bend left] coordinate (mp) (-2.5,2.5,-1)
to[bend right] (-2.5,-2.5,-1)
to[bend right] coordinate (mm) (2.5,-2.5,-1)
-- cycle;


\shade[opacity=0.4, draw=black] (2.5,-2.5,0) -- (2.5,2.5,0) -- (-2.5,2.5,0) -- (-2.5,-2.5,0) -- cycle;

\draw[-stealth, thick] (-0.5,0.08,0) -- coordinate[pos=.3] (f) (1.0,1.0,1);
\draw[-stealth] (-0.5,0.08,0) -- coordinate[pos=.3] (f) (1.0,1.0,0);
\draw[dotted] (1.0,1.0,1) -- (1.0,1.0,0);
\draw[solid]  (1.0,1.0,0) -- (1.0,1.0,0.25) -- (0.75,0.83,0.25) -- (0.75,0.83,0);
\filldraw (-0.5,0.08,0) circle (1pt);


\node[right, scale=0.8] at (-0.1,-0.5,0) {{$\ansatz(\param(t), \cdot)$}};
\node[right, scale=0.8] at (1.7,-0.3,0.2) {{$\partial_t\ansatz(\param(t), \cdot)$}};
\node[above, scale=0.8] at (1.0,1.0,1.0) {{$\vecf(\cdot, \ansatz(\param(t),\cdot))$}};
\draw (-1.5,2.5,0) node[above, rotate=-4] {\small{tangent space at $\ansatz(\param(t),\cdot)$}}; 
\draw (0,-4.1,0) node[right] {$\mathcal{M}$};

\end{tikzpicture}
	\caption{Neural Galerkin schemes are based on the Dirac--Frenkel variational principle \cite{Dirac1930Note,Frenkel1934Wave}, which defines the time derivative $\dot\param(t)$ of the parameter $\param(t)$ to be the orthogonal projection of the right-hand side function $\vecf( \cdot, \ansatz(\param(t), \cdot))$ onto the tangent space of the manifold $\manifold$ at the current solution field $\ansatz(\param(t), \cdot)$. The tangent space is spanned by the component functions of the gradient $\nabla_{\param}\ansatz(\param(t), \cdot)$.}
	\label{fig:manifold_illustration_plain}
\end{SCfigure}
Consider now a parametrization $\ansatz: \Theta \times \spacedomain \to \mathbb{R}^{\outputdim}$ of a solution field $\vecu$, which may depend nonlinearly on a time-dependent parameter vector $\param: \timedomain \rightarrow \Theta\subseteq\R^\paramdim$ of dimension $\paramdim$.
For example, the function $\ansatz$ can be a deep network, where the components of the parameter $\param(t) \in \Theta$ correspond to the weights and biases. We only consider parametrizations that are continuously differentiable in the parameter.
Plugging $\ansatz$ into~\eqref{eq:general_pde} leads to the residual function
\begin{align}
	\label{eq:ng_residual}
	\vecr_t(\param(t), \dot \param(t), \vecx) = \nablaansatz^\T \dot \param(t) - \vecf(\vecx, \ansatz(\param(t),\cdot))\,.
\end{align}
The time derivative $\dot \param(t)$ is then determined by following the Dirac--Frenkel variational principle~\cite{Dirac1930Note,Frenkel1934Wave,Kramer1981,Lubich2008} such that
\begin{align}
	\label{eq:ng_inner_product}
	\ip{
		\partialparami{i}\ansatz(\param(t),\; \cdot),\; \vecr_t(\param(t), \dot \param(t), \cdot)
	}_\nu
	=0\,, \qquad i = 1, \dots, \paramdim\,,
\end{align}
where $\partialparami{i}\ansatz$ is the $i$-th component function of the gradient $\nabla_{\param} \ansatz$ of $\ansatz$ with respect to the parameter $\param$.
A note on the history of the Dirac--Frenkel variational principle can be found in \cite[Section~3.8]{LasserL2020Computing}.
The solution $\dot\param(t)$ can be interpreted as determining an orthogonal projection $\nabla_{\param} \ansatz(\param(t), \cdot)^\T \dot\param(t)$ of the right-hand side function $\vecf(\cdot, \ansatz(\param(t),\cdot))$ evaluated at $\ansatz(\param(t),\cdot)$ onto the tangent space $\mathcal{T}_{\ansatz(\param(t), \cdot)}\manifold$ at the point $\ansatz(\param(t), \cdot)$ of the parametrization manifold
\begin{align}
	\label{eq:parametrization_manifold}
	\manifold = \left\{\ansatz(\veceta, \cdot) \,|\, \veceta \in \Theta \subseteq \R^{\paramdim}\right\} \subseteq \mathcal{U}\,,
\end{align}
which is illustrated in Figure~\ref{fig:manifold_illustration_plain}.

Conditions~\eqref{eq:ng_inner_product} can be rewritten in matrix form as
\begin{align}
	\label{eq:ng_ode}
	\vecM(\param(t)) \dot \param(t) = \vecF(\param(t)),
\end{align}
with the matrix $\vecM(\param(t))$ and vector $\vecF(\param(t))$ having the following components
\begin{subequations}
	\label{eq:ng_ode_parts}
	\begin{align}
		\vecM_{ij}(\param(t)) & = \ip{\partialparami{i}\ansatz(\param(t), \cdot),\; \partialparami{j}\ansatz(\param(t), \cdot)}_\nu, \quad i,j = 1, \dots, \paramdim\,, \\
		\vecF_i(\param(t))    & = \ip{\partialparami{i}\ansatz(\param(t), \cdot),\; \vecf(\cdot, \ansatz(\param(t), \cdot ))}_\nu, \quad i = 1, \dots,\paramdim\,.
	\end{align}
\end{subequations}
Following \cite{BrunaPV2022Neural}, we refer to \eqref{eq:ng_ode} as the Neural Galerkin system because the parametrizations that we use in the following are all based on neural networks and \eqref{eq:ng_inner_product} can be interpreted as a Galerkin projection with the component functions of $\nabla_{\param}\ansatz$ as test functions.
If the matrix elements \eqref{eq:ng_ode_parts} are not analytically available, then they can be numerically estimated via either quadrature or Monte Carlo methods.
We denote the sample-based estimates of $\vecM(\param(t))$ and vector $\vecF(\param(t))$ at the sampling points $\vecx_1, \dots, \vecx_{\nsamples} \in \spacedomain$ as
\begin{align}
	\hat\vecM(\param(t)) = \frac{1}{\nsamples}\sum\limits_{i=1}^{\nsamples} \nablaansatzi \nablaansatzi^\T, \label{eq:Prelim:NG:SampledM} \\
	\hat{\vecval{F}}(\param(t)) = \frac{1}{\nsamples}\sum\limits_{i=1}^{\nsamples} \nablaansatzi \vecval{f}(\cdot, \ansatz(\param(t),\cdot))(\vecx_i)\,,\label{eq:Prelim:NG:SampledF}
\end{align}
which give rise to the sampled Neural Galerkin system
\begin{equation}\label{eq:Prelim:SampledNGODE}
	\hat\vecM(\param(t)) \dot \param(t) = \hat{\vecF}(\param(t))\,.
\end{equation}
There are several approaches for obtaining efficient Monte Carlo estimates \cite{BrunaPV2022Neural,WenVP2023Coupling}, which goes beyond the scope of the present work.
In the following, we assume that the system of ordinary differential equations (ODEs) given by the Neural Galerkin system \eqref{eq:ng_ode} and the corresponding sampled system \eqref{eq:Prelim:SampledNGODE} have a solution on the whole time interval $\timedomain$.

\subsection{Problem formulation}
\label{sec:problem_statement}
Conserving quantities in  solutions obtained with Neural Galerkin schemes based on nonlinear parametrizations  leads to two challenges.
First, a quantity $\cq_i$ is not necessarily conserved by a Neural Galerkin solution $\ansatz(\param(t), \cdot)$ with parameter $\param(t)$ satisfying \eqref{eq:ng_ode}, even in continuous time and without sampling. The reason is that the Galerkin projection formulated in \eqref{eq:ng_inner_product} only seeks to set the residual orthogonal to a tangent space of $\manifold$ but ignores any additional constraints given by the quantities  $\cq_1, \dots, \cq_{\ncq}$.

Second, only adding constraints to the continuous-time Neural Galerkin equations \eqref{eq:ng_ode} or their sampled counterparts \eqref{eq:Prelim:SampledNGODE} is insufficient to conserve quantities because of the nonlinear parametrization of the solution field.
While linear quantities can be conserved by implicit and explicit Runge-Kutta integrators~\cite{Shampine1986Conservation} and quadratic quantities can be conserved by Runge-Kutta integrators if the coefficients satisfy conditions described in~\cite{Cooper1987Stability}, arbitrary higher-order polynomial or nonlinear quantities are not conserved by any Runge-Kutta integrators in general; see also the discussion in~\cite[Chapter IV]{HairerWL2006Geometric}.
This is important in the context of Neural Galerkin schemes and related methods because  the parameter $\param(t)$ can enter nonlinearly in the parametrization $\ansatz(\param(t), \cdot)$. This means that conserved quantities that are linear in the solution field $\vecu$ of the PDE formulation~\eqref{eq:general_pde} (e.g., the integral in Burgers' equation) depend through the parametrization $\ansatz(\param(t), \cdot)$ nonlinearly on the parameter $\param(t)$, which is the state of the Neural Galerkin system \eqref{eq:ng_ode} and \eqref{eq:Prelim:SampledNGODE}.
Thus, even if the continuous-time Neural Galerkin projection is constrained so that $\dot \param(t)$ conserves the quantities $\cq_1, \dots, \cq_{\ncq}$, the conservation is lost after discretizing in time with common integrators such as Runge-Kutta integrators because they typically only conserve linear quantities. In fact, there is no Runge-Kutta method that can conserve all polynomial quantities of degree greater than two~\cite[Theorem IV.3.3]{HairerWL2006Geometric}.

\section{Conserving quantities in continuous-time Neural Galerkin schemes}\label{sec:conserved_ng}
We pursue two options for conserving quantities in continuous time. First, in Section~\ref{sec:StrucConstraint}, we conserve quantities via constraints by building on previous work and introducing constrained Neural Galerkin schemes that conserve quantities in continuous time for generic nonlinear parametrizations. Second, in Section~\ref{sec:StrucHam}, we enforce conservation via structure in the nonlinear parametrization. In particular, we construct nonlinear parametrizations so that weighted Neural Galerkin schemes preserve Hamiltonians without the need of explicitly adding constraints.

\subsection{Adding constraints to Neural Galerkin schemes for conserving quantities in continuous time}\label{sec:StrucConstraint}
\begin{figure}[t]
	\centering
	\begin{tabular}{cc}
		\begin{tikzpicture}[scale=0.93]
\begin{scope}[scale=0.9]

\tikzset{x={(175:1cm)},y={(55:.7cm)},z={(90:1cm)}}

\shade[left color=blue!4,right color=blue!60,looseness=.5, draw=black]  (2.5,-2.5,-1) 
to[bend left] (2.5,2.5,-1)
to[bend left] coordinate (mp) (-2.5,2.5,-1)
to[bend right] (-2.5,-2.5,-1)
to[bend right] coordinate (mm) (2.5,-2.5,-1)
-- cycle;

\draw[red, thick] (-2.5,0,-0.55) .. controls (-1.5,0,0) and (0.5,1.3,0) .. (2.5,-1.7, -0.62);

\shade[opacity=0.4, draw=black] (2.5,-2.5,0) -- (2.5,2.5,0) -- (-2.5,2.5,0) -- (-2.5,-2.5,0) -- cycle;

\draw[-stealth, thick] (-0.5,0.08,0) -- coordinate[pos=.3] (f) (1.0,1.0,1);
\draw (-0.5,0.08,0) -- coordinate[pos=.3] (f) (1.0,1.0,0);
\draw[dotted] (1.0,1.0,1) -- (1.0,1.0,0);
\draw[solid]  (1.0,1.0,0) -- (1.0,1.0,0.25) -- (0.75,0.83,0.25) -- (0.75,0.83,0);
\draw[dotted] (1.0,1.0,0) -- (1.0,0.08,0);
\draw[solid]  (1.0,0.08,0) -- (1.0,0.33,0.0) -- (0.75,0.33,0.0) -- (0.75,0.08,0) -- cycle;
\draw[red] (-2.5,0.08,0) -- (2.5,0.08,0);
\draw[-stealth, dashed]  (-0.5, 0.08, 0.0) -- (1.0, 0.08, 0.0);
\filldraw (-0.5,0.08,0) circle (1pt);


\node[above, scale=0.8] at (1.0,1.0,1.0) {$\vecf(\cdot, \ansatz(\param(t),\cdot))$};
\draw (0,-4.1,0) node[right] {$\mathcal{M}$};
\draw (-3.0,-1.7,0) node[right] {\textcolor{red}{$\mathcal{M}_{I}$}};

\end{scope}
\end{tikzpicture}                       &
		\begin{tikzpicture}[scale=0.93]
\begin{scope}[scale=0.95]

\tikzset{x={(175:1cm)},y={(55:.7cm)},z={(90:1cm)}}

\shade[left color=blue!4,right color=blue!60,looseness=.5, draw=black]  (2.5,-2.5,-1) 
to[bend left] (2.5,2.5,-1)
to[bend left] coordinate (mp) (-2.5,2.5,-1)
to[bend right] (-2.5,-2.5,-1)
to[bend right] coordinate (mm) (2.5,-2.5,-1)
-- cycle;

\draw[red, thick] (-2.5,0,-0.55) .. controls (-1.5,0,0) and (0.5,1.3,0) .. (2.5,-1.5, -0.62);


\draw[-stealth, solid, looseness=0.7, thick]  (-0.5, 0.08, 0.0) .. controls (0.5, 0.00, 0.1) ..(1.5, 0.12, -0.2);
\filldraw (-0.5,0.08,0) circle (1pt);
\filldraw (1.5,0.12,-0.2) circle (1pt);
\filldraw (1.1,-0.77,0) circle (1pt);
\draw[dotted]  (1.5, 0.12, -0.2) -- (1.1,-0.77,0);


\draw (0,-4.1,0) node[right] {$\mathcal{M}$};

\node[above, scale=0.7] at (-0.75,-0.08,0.1) {$\ansatz(\param_k, \cdot)$};
\node[above, scale=0.7, rotate=8.5] at (1.1,0.12,-0.1) {$\ansatz(\tilde{\param}_{k + 1}, \cdot)$};
\node[below, scale=0.7] at (0.5,-0.77,-0.05) {$\ansatz(\param_{k + 1}, \cdot)$};
\draw (-3.0,-1.7,0) node[right] {\textcolor{red}{$\mathcal{M}_{I}$}};
\end{scope}

\end{tikzpicture}                                                                   \\
		\vspace{0.1cm}                                                 &                                                  \\
		\scriptsize{(a) adding constraints to Neural Galerkin schemes} & \scriptsize{(b) Neural Galerkin with embeddings}
	\end{tabular}
	\caption{Only adding constraints to the continuous-time Neural Galerkin
		system is insufficient to conserve
		quantities because the constraint can be violated in discrete time due to the nonlinear parametrization of the solution field. In contrast, the proposed Neural Galerkin scheme combines constraints with nonlinear projections to obtain embeddings onto manifolds of functions that conserve quantities in discrete time.}
	\label{fig:manifold_illustration}
\end{figure}
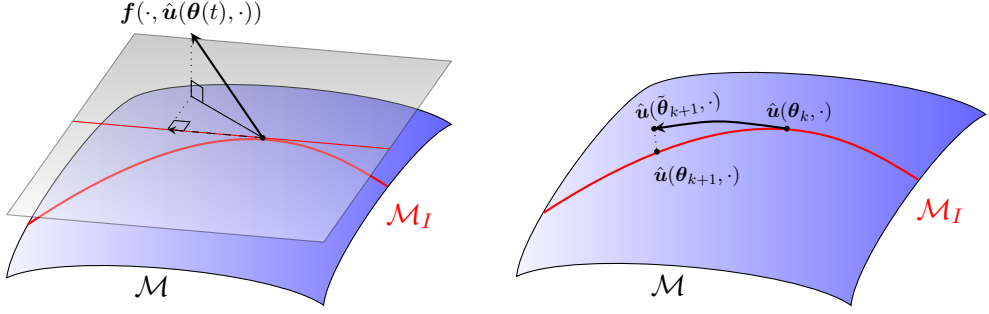
The Neural Galerkin scheme is based on projecting the right-hand side onto tangent spaces of the parametrization manifold; see Section~\ref{sec:NG}.
If the image of the right-hand side function is a subset of the tangent space, then the residual is zero and thus all quantities that are conserved by the PDE solution field $\vecu$ are also conserved by the continuous-time Neural Galerkin solution $\ansatz(\param(t), \cdot)$. This relation between residual and conserving quantities has been, for instance, used in the context of the Schr\"odinger equation in \cite[Sec.~3]{LasserL2020Computing} to establish energy and norm conservation when using a Gaussian wave packet for the parametrization.
We now consider more general cases where the residual is not necessarily zero and add constraints to the time-continuous Neural Galerkin system \eqref{eq:ng_ode} and its sampled counterpart \eqref{eq:Prelim:SampledNGODE} to conserve quantities in continuous time; similar to other methods based on nonlinear parametrizations \cite{LeeC2021Deep,AndersonF2022Evolution}.

\subsubsection{Constrained manifolds}
Recall the interpretation that $\dot\param(t)$ corresponds to an orthogonal projection of the right-hand side function $\vecf$ onto a tangent space of $\manifold$.
By adding a constraint to the Neural Galerkin system, we restrict the manifold $\manifold$ of parametrized functions $\ansatz(\veceta, \cdot)$ defined in \eqref{eq:parametrization_manifold} to the manifold $\manifoldC$ of functions that conserve quantities
\begin{align}
	\label{eq:manifold_constraint_constant}
	\cq_i(\ansatz(\veceta, \cdot)) = \manifoldconstant_i\,,\qquad i = 1, \dots, \ncq\,,
\end{align}
for a given vector of constants $\vecval{\manifoldconstant} = [\manifoldconstant_1, \dots, \manifoldconstant_\ncq]$.
Thus, we obtain $\manifoldC$ as
\begin{align}
	\label{eq:ng_manifold_def}
	\manifoldC & = \{ \ansatz(\veceta, \cdot) : \veceta \in \Theta \text{ and } \cq_i(\ansatz(\veceta, \cdot)) = \manifoldconstant_i \text{ for } i =1, \dots, \ncq\}.
\end{align}

\subsubsection{Constrained Neural Galerkin schemes}
We now apply the Dirac--Frenkel variational principle with respect to the constrained manifold $\manifoldC$ instead of $\manifold$ to determine $\dot\param(t)$; see Figure~\ref{fig:manifold_illustration}a.
The projection onto a tangent space of $\manifoldC$ is analogous to the projection onto a tangent space of $\manifold$ and leads to the constrained Neural Galerkin system
\begin{align}
	\label{eq:modified_ng_ode}
	\begin{bmatrix}
		\vecM(\param(t))    & \vecg(\param(t)) \\
		\vecg(\param(t))^\T & 0
	\end{bmatrix}
	\begin{bmatrix}
		\dot \param(t) \\
		\veclambda(t)
	\end{bmatrix}
	=
	\begin{bmatrix}
		\vecF(\param(t)) \\
		0
	\end{bmatrix},
\end{align}
where
$\vecg(\param(t))$ is
\begin{equation}\label{eq:CNG:GVecCont}
	\vecg(\param(t)) =\begin{bmatrix}
		\nablaparam \cq_1(\ansatz(\param(t), \cdot)), \dots, \nablaparam \cq_{\ncq}(\ansatz(\param(t), \cdot))
	\end{bmatrix}\,.
\end{equation}
System \eqref{eq:modified_ng_ode} is analogous to the system introduced in \cite{AndersonF2022Evolution}. The vector $\veclambda(t)$ contains the Lagrange multipliers at time $t$.
The system \eqref{eq:modified_ng_ode}
has a unique solution $(\dot \param(t),\veclambda(t))$ if $\vecM(\param(t))$ is regular and $\vecg(\param(t))$ has full column rank~\cite[Chapter 5]{Bjorck1996Numerical}.
Moreover, by verifying that the requirements of \cite[Thm.~4.13]{KunkelM2006Differential-Algebraic} are satisfied, the initial value problem associated with the system of differential-algebraic equations \eqref{eq:modified_ng_ode} has locally a unique solution, provided that $\vecM$ is pointwise invertible, $\vecg$ has pointwise full column rank, and $\vecM,\vecg,\vecF$ are continuously differentiable.
If the constant in the definition \eqref{eq:ng_manifold_def} of the manifold $\manifoldC$ is set to $\vecval{\manifoldconstant} = \begin{bmatrix} \cq_1(\ansatz(\param(0), \cdot), \dots, \cq_\ncq(\ansatz(\param(0), \cdot) \end{bmatrix}$, then a continuous-time Neural Galerkin solution $\ansatz(\param(t), \cdot)$ with $\param(t)$ satisfying the constrained Neural Galerkin system \eqref{eq:modified_ng_ode} conserves the quantities $\cq_1, \dots, \cq_{\ncq}$.

\subsubsection{Sampled constrained Neural Galerkin schemes}
Analogously, we can introduce the manifold $\hat\manifoldC$ that is based on the sampled quantities \eqref{eq:conserved_quantity_monte_carlo},
\begin{equation}\label{eq:NonlinProj:SampledME}
	\hat\manifoldC =  \{\ansatz(\veceta, \cdot) : \veceta \in \Theta \text{ and } \hat\cq_i(\ansatz(\veceta, \cdot)) = \manifoldconstant_i \text{ for } i =1, \dots, \ncq\}\,,
\end{equation}
and derive the sampled constrained Neural Galerkin system
as \begin{equation}\label{eq:modified_ng_ode_sampled}
	\begin{bmatrix}
		\hat{\vecM}(\param(t))    & \hat{\vecg}(\param(t)) \\
		\hat{\vecg}(\param(t))^\T & 0
	\end{bmatrix}
	\begin{bmatrix}
		\dot \param(t) \\
		\veclambda(t)
	\end{bmatrix}
	=
	\begin{bmatrix}
		\hat{\vecF}(\param(t)) \\
		0
	\end{bmatrix},
\end{equation}
where $\hat{\vecM}(\param(t))$ and $\hat{\vecF}(\param(t))$ are the sampled $\vecM(\param(t))$ and $\vecF(\param(t))$, respectively. In $\hat{\vecg}(\param(t))$, the sampled quantities \eqref{eq:conserved_quantity_monte_carlo} are used as
\begin{equation}\label{eq:TimeCont:SampledGVec}
	\hat{\vecg}(\param(t)) =\begin{bmatrix}
		\nablaparam \hat{\cq}_1(\ansatz(\param(t), \cdot)), \dots, \nablaparam \hat{\cq}_{\ncq}(\ansatz(\param(t), \cdot))
	\end{bmatrix}\,,
\end{equation}
so that solutions of the sampled constrained system \eqref{eq:modified_ng_ode_sampled} conserve the sampled quantities \eqref{eq:conserved_quantity_monte_carlo}.

\subsection{Structured nonlinear parametrizations and weighted Neural Galerkin schemes to preserve Hamiltonians in continuous time}\label{sec:StrucHam}
We now derive specific nonlinear parametrizations and weighted schemes that preserve Hamiltonians without having to resort to using constraints.

\subsubsection{Separable nonlinear parametrizations}
To preserve Hamiltonians in time-continuous Neural Galerkin solutions, we consider separable parametrizations that are of the form
\begin{equation}
	\label{eq:separable_ansatz}
	\ansatz(\param(t), \vecx) = \sum\nolimits_{i=1}^{\nPhi} \vecbeta_i(t) \phi_i(\vecx, \vecalpha_i(t))
\end{equation}
with $\vecbeta_i: \timedomain \rightarrow \R^{\outputdim}$, $\vecalpha_i: \timedomain \rightarrow \R^{\alphadim_i}$, $\phi_i: \spacedomain \times \R^{\alphadim_i} \rightarrow \R$. The parameters can be combined into
\begin{equation}\label{eq:ContTime:SepStructParam}
	\param(t) = [\vecbeta_1(t)^\T, \dots, \vecbeta_{\nPhi}(t)^\T, \vecalpha_1(t)^\T, \dots, \vecalpha_{\nPhi}(t)^\T]^\T \in \Theta\,.
\end{equation}
We have $\paramdim = \nPhi\outputdim+\sum\nolimits_{i=1}^{\nPhi}\alphadim_i$ and $\Theta = \R^{\paramdim}$.
A parametrization of the form \eqref{eq:separable_ansatz} is separable because $\param$ can be separated into the components $\vecalpha$ and $\vecbeta$, where $\vecbeta$ enters the parametrization linearly.
A similar separable parametrization has been used in the context of nonlinear model reduction for finite-dimensional port-Hamiltonian systems in \cite{Schulze2023Structure-Preserving}.
In the context of deep-network parametrizations, an architecture \eqref{eq:separable_ansatz} is obtained whenever the last (output) layer of the network is linear without bias.
In the following, it will be convenient to introduce the matrix function $\vecV\colon \spacedomain\times \R^{\sum\nolimits_{i=1}^{\nPhi}\alphadim_i} \to \R^{\outputdim\times \nPhi\outputdim}$
via $\vecV(\vecx,\vecalpha) =
	[
	\phi_1(\vecx, \vecalpha_1), \cdots, \phi_{\nPhi}(\vecx, \vecalpha_{\nPhi})
	]
	\otimes \boldsymbol{1}_{\outputdim}$,
where $\otimes$ denotes the Kronecker product and $\boldsymbol{1}_{\outputdim}$ the $\outputdim \times \outputdim$ identity matrix.

\subsubsection{Hamiltonians with factorizable structure}\label{sec:special_hamiltonians}
We call a Hamiltonian $\hamiltonian$ factorizable if there exist a continuously differentiable function $h\colon \R^{\outputdim}\to\R$ and a point-wise symmetric and positive definite matrix function $\vecQ\colon \R^{\outputdim}\to \R^{\outputdim\times\outputdim}$ that satisfy
\begin{align}
	\label{eq:assumption_hamiltonian_depends_only_on_state_values}
	\hamiltonian(\vecv)                                               & = \int_{\spacedomain} h(\vecv(\vecx)) \dd \nu(\vecx)\,,\qquad \text{ and }                            \\
	\label{eq:assumptionForVariationalDerivative}
	\left(\frac{\delta\hamiltonian}{\delta\vecu}(\vecv)\right)(\vecx) & = \vecQ(\vecv(\vecx))\vecv(\vecx)\quad\text{for all }(\vecv,\vecx)\in\mathcal{U}\times\spacedomain\,,
\end{align}
where $\delta H/\delta \vecu$ denotes the variational derivative of $H$.
Hamiltonians are factorizable if they correspond to a squared norm of the state, which is, for instance, the case for the Burgers' and wave equation examples in \eqref{eq:burgers_equation_hamiltonian_structure}--\eqref{eq:wave_equation_hamiltonian_structure}.
However, the Hamiltonian considered for the shallow water equations in the later Section~\ref{sec:shallow_water} is not factorizable in the sense of \eqref{eq:assumption_hamiltonian_depends_only_on_state_values}--\eqref{eq:assumptionForVariationalDerivative} because it involves the gradient of one of the state variables.

\subsubsection{Weighted Neural Galerkin schemes that preserve Hamiltonians in continuous time}
We now introduce a Neural Galerkin scheme that conserves Hamiltonians that satisfy \eqref{eq:assumptionForVariationalDerivative} with a matrix function $\vecQ$.
To this end, we follow the Dirac--Frenkel approach as in \eqref{eq:ng_inner_product} and perform the projection with respect to the weighted inner product given by the function $\vecQ$ of \eqref{eq:assumptionForVariationalDerivative}; see also \cite{Schulze2023Structure-Preserving} for a similar approach in the context of model reduction.
Consider a separable nonlinear parametrization \eqref{eq:separable_ansatz} with the parameter $\param(t)$ given in \eqref{eq:ContTime:SepStructParam} and define $\vecM_{\vecQ}\colon \R^{\paramdim}\to\R^{\paramdim\times\paramdim}$ and $\vecF_{\vecQ}\colon \R^{\paramdim}\to\R^{\paramdim}$ as
\begin{equation}		      \label{eq:structured_neural_galerkin_coefficients}
	\vecM_{\vecQ}(\param(t)) =
	\begin{bmatrix}
		\vecM^{(11)}(\param(t))    & \vecM^{(12)}(\param(t)) \\
		\vecM^{(12)}(\param(t))^\T & \vecM^{(22)}(\param(t))
	\end{bmatrix}
	,\qquad \vecF_{\vecQ}(\param(t)) =
	\begin{bmatrix}
		\vecF^{(1)}(\param(t)) \\
		\vecF^{(2)}(\param(t))
	\end{bmatrix}\,,
\end{equation}
with the blocks defined as
\begin{equation}
	\begin{aligned}
		\vecM^{(11)}_{ij}(\param(t)) =    & \ip{\vecQ(\ansatz(\param(t),\cdot))\partial_{\vecalpha}\vecV(\cdot,\vecalpha(t))(\vece_i)\vecbeta(t),\; \partial_{\vecalpha}\vecV(\cdot,\vecalpha(t))(\vece_j)\vecbeta(t)}_{\nu}\,, \\
		\vecM^{(12)}_{i\ell}(\param(t)) = & \ip{\vecQ(\ansatz(\param(t),\cdot))\partial_{\vecalpha}\vecV(\cdot,\vecalpha(t))(\vece_i)\vecbeta(t),\; \vecV(\cdot,\vecalpha(t))\vece_\ell}_{\nu}\,,                               \\
		\vecM^{(22)}_{k\ell}(\param(t)) = & \ip{ \vecQ(\ansatz(\param(t),\cdot))\vecV(\cdot,\vecalpha(t))\vece_k,\; \vecV(\cdot,\vecalpha(t))\vece_\ell}_{\nu}\,,                                                               \\
		\vecF^{(1)}_i(\param(t)) =        & \left\langle \vecQ(\ansatz(\param(t),\cdot))\partial_{\vecalpha}\vecV(\cdot,\vecalpha(t))(\vece_i)\vecbeta(t),\right.                                                               \\
		                                  & \qquad\qquad\qquad\qquad\left.\;\Joperator(\ansatz(\param(t),\cdot))\vecQ(\ansatz(\param(t),\cdot))\ansatz(\param(t),\cdot)\right\rangle_{\nu}\,,                                   \\
		\vecF^{(2)}_k(\param(t)) =        & \ip{ \vecQ(\ansatz(\param(t),\cdot))\vecV(\cdot,\vecalpha(t))\vece_k,\;\Joperator(\ansatz(\param(t),\cdot))\vecQ(\ansatz(\param(t),\cdot))\ansatz(\param(t),\cdot)}_{\nu}\,,
	\end{aligned}
\end{equation}
for $i,j=1,\ldots,\sum_{s=1}^{\nPhi}\alphadim_s$, $k,\ell=1,\ldots,\nPhi\outputdim$, and the $i$th canonical unit vector $\vece_i$.
The matrix functions $\vecM_{\vecQ}$ and $\vecF_{\vecQ}$ lead to the weighted time-continuous Neural Galerkin system
\begin{equation}
	\label{eq:NonautonomousSemiDiscretizedSystem}
	\vecM_{\vecQ}(\param(t))\dot{\param}(t) = \vecF_{\vecQ}(\param(t))\,.
\end{equation}
The following proposition states that $\hamiltonian$ is indeed a conserved quantity of \eqref{eq:NonautonomousSemiDiscretizedSystem}, where we overload the notation to use $\hamiltonian(\param)$ as short-hand notation for $\hamiltonian(\ansatz(\param,\cdot))$.
\begin{proposition}
	\label{thm:structure-preserving_neural_galerkin}
	Consider a Hamiltonian PDE of the form \eqref{eq:general_pde} with right-hand side satisfying \eqref{eq:structuredRHS}--\eqref{eq:skew-adjointness_J} and let the Hamiltonian be factorizable. Then, any solution of the corresponding time-continuous weighted Neural Galerkin system 
	\eqref{eq:NonautonomousSemiDiscretizedSystem} based on the separable parametrization \eqref{eq:separable_ansatz} preserves the Hamiltonian $\hamiltonian$ in the sense of
	\begin{equation*}
		\frac{\dd\hamiltonian(\param(t))}{\dd t}(t) = 0\quad \text{for all }t\in\timedomain.
	\end{equation*}
\end{proposition}
\begin{proof}
	Straightforward calculations lead to the observations that the Hamiltonian of the time-continuous Neural Galerkin system satisfies
	\begin{equation}
		\label{eq:gradient_of_hamiltonian_of_discretized_system}
		\nabla\hamiltonian(\veceta) = \vecM_{\vecQ}(\veceta)^\T
		\begin{bmatrix}
			0 & 0                                \\
			0 & \boldsymbol{1}_{\nPhi\outputdim}
		\end{bmatrix}
		\veceta
	\end{equation}
	and that the right-hand side of \eqref{eq:NonautonomousSemiDiscretizedSystem} satisfies
	\begin{equation}
		\label{eq:factorization_of_structured_rhs}
		\vecF(\veceta) = \vecJ_{\vecQ}(\veceta)
		\begin{bmatrix}
			0 & 0                                \\
			0 & \boldsymbol{1}_{\nPhi\outputdim}
		\end{bmatrix}
		\veceta
	\end{equation}
	for all $\veceta\in\R^{\paramdim}$, where $\vecJ_{\vecQ}\colon \R^{\paramdim}\to\R^{\paramdim\times\paramdim}$ is defined via
	\begin{equation*}
		\vecJ_{\vecQ}(\param(t)) :=
		\begin{bmatrix}
			0                           & \vecJ^{(12)}(\param(t)) \\
			-\vecJ^{(12)}(\param(t))^\T & \vecJ^{(22)}(\param(t))
		\end{bmatrix}\,,
	\end{equation*}
	with the blocks
	\begin{align}
		\phantom{\vecJ_{k,\ell}^{(22)}(\param(t))}
		 & \begin{aligned}
			   \mathllap{\vecJ^{(12)}_{i,\ell}(\param(t))} = & \left\langle \vecQ(\ansatz(\param(t),\cdot))\partial_{\vecalpha}\vecV(\cdot,\vecalpha(t))(\vece_i)\vecbeta(t),\right.                       \\
			                                                 & \qquad\quad\left.\Joperator(\ansatz(\param(t),\cdot))\vecQ(\ansatz(\param(t),\cdot))\vecV(\cdot,\vecalpha(t))\vece_\ell\right\rangle_{\nu},
		   \end{aligned} \\
		 & \begin{aligned}
			   \mathllap{\vecJ_{k,\ell}^{(22)}(\param(t))} = & \left\langle \vecQ(\ansatz(\param(t),\cdot))\vecV(\cdot,\vecalpha(t))\vece_k,\right.                                                        \\
			                                                 & \qquad\quad\left.\Joperator(\ansatz(\param(t),\cdot))\vecQ(\ansatz(\param(t),\cdot))\vecV(\cdot,\vecalpha(t))\vece_\ell\right\rangle_{\nu},
		   \end{aligned}\label{eq:StrucHam:J22}
	\end{align}
	for $i=1,\ldots,\sum_{s=1}^{\nPhi}\alphadim_s$ and $k,\ell=1,\ldots,\nPhi\outputdim$.
	By exploiting the pointwise skew-adjointness of $\Joperator$, we conclude that $\vecJ^{(22)}$ and $\vecJ_{\vecQ}$ are pointwise skew-symmetric.
	The conservation of the Hamiltonian follows then from \cite{MehrmannM2019Structure-Preserving} by observing that \eqref{eq:gradient_of_hamiltonian_of_discretized_system} and \eqref{eq:factorization_of_structured_rhs} imply that \eqref{eq:structured_neural_galerkin_coefficients}--\eqref{eq:NonautonomousSemiDiscretizedSystem} has a port-Hamiltonian structure without dissipation or input/output ports.
\end{proof}

\subsubsection{Weighted Neural Galerkin schemes with Monte Carlo approximations}
We now show that the sampled quantities are conserved with weighted Neural Galerkin schemes in continuous time even if the integrals occurring in $\vecM_{\vecQ}$ and $\vecF_{\vecQ}$ are replaced by corresponding Monte Carlo estimates.
Replacing $\vecM_{\vecQ}$ with its sampled counterpart $\hat{\vecM}_{\vecQ}$ analogous to $\hat{\vecM}$ defined in \eqref{eq:Prelim:NG:SampledM} poses no problems in terms of structure preservation. In contrast, the right-hand side term $\vecF_{\vecQ}$ needs more careful treatment because applying Monte Carlo directly to $\vecF_{\vecQ}$ can destroy the property \eqref{eq:factorization_of_structured_rhs} that is used in the proof of Proposition~\ref{thm:structure-preserving_neural_galerkin}.
To avoid loosing property \eqref{eq:factorization_of_structured_rhs} in the sampled $\vecF_{\vecQ}$, we first use \eqref{eq:factorization_of_structured_rhs} and the skew symmetry of $\vecJ_{22}$ defined in \eqref{eq:StrucHam:J22} and write $\vecF_{\vecQ}$ as
\begin{align*}
	\vecF_{\vecQ}(\param(t)) & =
	\begin{bmatrix}
		0                           & \vecJ^{(12)}(\param(t))                                     \\
		-\vecJ^{(12)}(\param(t))^\T & \frac12(\vecJ^{(22)}(\param(t))-\vecJ^{(22)}(\param(t))^\T)
	\end{bmatrix}
	\begin{bmatrix}
		0 & 0                                \\
		0 & \boldsymbol{1}_{\nPhi\outputdim}
	\end{bmatrix}
	\param(t)                    \\
	                         & =
	\begin{bmatrix}
		\vecJ^{(12)}(\param(t))\vecbeta(t) \\
		\frac12\vecJ^{(22)}(\param(t))\vecbeta(t)-\frac12\vecJ^{(22)}(\param(t))^\T\vecbeta(t)
	\end{bmatrix}
	.
\end{align*}
Instead of approximating $\vecF_{\vecQ}(\param(t))$ directly via Monte Carlo estimates, we propose to approximate $\vecJ^{(12)}(\param(t))\vecbeta(t)$ as well as $\vecJ^{(22)}(\param(t))\vecbeta(t)$ and $\vecJ^{(22)}(\param(t))^\T\vecbeta(t)$ separately.
The approximations can then be used to assemble the approximation $\hat{\vecF}_{\vecQ}$ of $\vecF_{\vecQ}$ that can be factorized analogous to \eqref{eq:factorization_of_structured_rhs}.
The sampled weighted Neural Galerkin system is
\begin{equation}
	\label{eq:structuredSampledSystem}
	\Msampled(\param(t))\dot\param(t) = \Jsampled(\param(t))
	\begin{bmatrix}
		0 & 0                                \\
		0 & \boldsymbol{1}_{\nPhi\outputdim}
	\end{bmatrix}
	\param(t)\,,
\end{equation}
where the Monte Carlo approximations of $\vecM_{\vecQ}$ and $\vecJ_{\vecQ}$ are denoted as $\Msampled,\Jsampled\colon \R^{\paramdim}\to\R^{\paramdim\times\paramdim}$, respectively, and defined as
\begin{equation}
	\label{eq:sampled_structured_neural_galerkin_coefficients}
	\begin{aligned}
		\Msampled(\param(t)) = &
		\begin{bmatrix}
			\Msampledblock{11}(\param(t))    & \Msampledblock{12}(\param(t)) \\
			\Msampledblock{12}(\param(t))^\T & \Msampledblock{22}(\param(t))
		\end{bmatrix}
		,                        \\
		\Jsampled(\param(t)) = &
		\begin{bmatrix}
			0                                 & \Jsampledblock{12}(\param(t)) \\
			-\Jsampledblock{12}(\param(t))^\T & \Jsampledblock{22}(\param(t))
		\end{bmatrix}
	\end{aligned}
\end{equation}
with the blocks
\begin{align*}
	\phantom{\Msampledblock{22}_{k\ell}(\param(t))}
	 & \begin{aligned}
		   \mathllap{\Msampledblock{11}_{ij}(\param(t))} = & \frac1{\nsamples}\sum_{s=1}^{\nsamples}\vecbeta(t)^\T\left(\partial_{\vecalpha}\vecV(\vecx_s,\vecalpha(t))\vece_i\right)^\T \vecQ(\ansatz(\param(t),\vecx_s))^\T \\
		                                                   & \qquad\quad\cdot\partial_{\vecalpha}\vecV(\vecx_s,\vecalpha(t))(\vece_j)\vecbeta(t),
	   \end{aligned}        \\
	 & \begin{aligned}
		   \mathllap{\Msampledblock{12}_{i\ell}(\param(t))} = & \frac1{\nsamples}\sum_{s=1}^{\nsamples}\vecbeta(t)^\T\left(\partial_{\vecalpha}\vecV(\vecx_s,\vecalpha(t))\vece_i\right)^\T \vecQ(\ansatz(\param(t),\vecx_s))^\T \\
		                                                      & \qquad\quad\cdot\vecV(\vecx_s,\vecalpha(t))\vece_\ell,
	   \end{aligned}     \\
	 & \begin{aligned}
		   \mathllap{\Msampledblock{22}_{k\ell}(\param(t))} = & \frac1{\nsamples}\sum_{s=1}^{\nsamples}\vece_k^\T \vecV(\vecx_s,\vecalpha(t))^\T \vecQ(\ansatz(\param(t),\vecx_s))^\T \vecV(\vecx_s,\vecalpha(t))\vece_\ell,
	   \end{aligned} \\
	 & \begin{aligned}
		   \mathllap{\Jsampledblock{12}_{i\ell}(\param(t))}  = & \frac1{\nsamples}\sum_{s=1}^{\nsamples}\vecbeta(t)^\T\left(\partial_{\vecalpha}\vecV(\vecx_s,\vecalpha(t))\vece_i\right)^\T \vecQ(\ansatz(\param(t),\vecx_s))^\T \\
		                                                       & \qquad\quad\cdot\left(\Joperator(\ansatz(\param(t),\cdot))\vecQ(\ansatz(\param(t),\cdot))\vecV(\cdot,\vecalpha(t))\vece_\ell\right)(\vecx_s)\,,
	   \end{aligned}                   \\
	 & \begin{aligned}
		   \mathllap{\Jsampledblock{22}_{k\ell}(\param(t))} = & \frac1{2\nsamples}\sum_{s=1}^{\nsamples}\vece_k^\T \vecV(\vecx_s,\vecalpha(t))^\T  \vecQ(\ansatz(\param(t),\vecx_s))^\T                      \\
		                                                      & \qquad\quad\cdot\left(\Joperator(\ansatz(\param(t),\cdot))\vecQ(\ansatz(\param(t),\cdot))\vecV(\cdot,\vecalpha(t))\vece_\ell\right)(\vecx_s) \\
		                                                      & -\frac1{2\nsamples}\sum_{s=1}^{\nsamples}\vece_\ell^\T \vecV(\vecx_s,\vecalpha(t))^\T  \vecQ(\ansatz(\param(t),\vecx_s))^\T                  \\
		                                                      & \qquad\quad\cdot\left(\Joperator(\ansatz(\param(t),\cdot))\vecQ(\ansatz(\param(t),\cdot))\vecV(\cdot,\vecalpha(t))\vece_k\right)(\vecx_s)\,,
	   \end{aligned}
\end{align*}
for $i,j=1,\ldots,\sum_{s=1}^{\nPhi}\alphadim_s$, $k,\ell=1,\ldots,\nPhi\outputdim$.

The following proposition states that the sampled Hamiltonian $\Hsampled$ is a conserved quantity of the sampled weighted Neural Galerkin system \eqref{eq:structuredSampledSystem}.
\begin{proposition}
	\label{thm:structure-preserving_neural_galerkin_with_sampling}
	Let the assumptions of Proposition~\ref{thm:structure-preserving_neural_galerkin} be satisfied and consider the  sampled weighted Neural Galerkin system \eqref{eq:structuredSampledSystem}.
	Moreover, let the associated Hamiltonian $\Hsampled$ be based on the same sampling points, i.e, $\nsamples=\nsamplesmani$ and $\vecx_i=\vecxi_i$ for $i=1,\ldots,\nsamples$.
	Then, any solution of \eqref{eq:structuredSampledSystem} satisfies
	\begin{equation*}
		\frac{\dd\Hsampled(\param(t))}{\dd t}(t) = 0\quad \text{for all }t\in\timedomain.
	\end{equation*}
\end{proposition}
\begin{proof}
	First, we note that the integral form \eqref{eq:assumption_hamiltonian_depends_only_on_state_values} of $\hamiltonian$ implies
	\begin{equation*}
		\frac{\delta\hamiltonian}{\delta \vecu}(\vecv) = \nabla h\circ \vecv = \vecQ(\vecv)\vecv
	\end{equation*}
	for all $\vecv\in\mathcal{U}$.
	Then, by straightforward calculations, we obtain that
	\begin{equation*}
		\nabla\Hsampled(\veceta) = \Msampled(\veceta)^\T
		\begin{bmatrix}
			0 & 0                                \\
			0 & \boldsymbol{1}_{\nPhi\outputdim}
		\end{bmatrix}
		\veceta
	\end{equation*}
	holds for all $\veceta\in\R^{\paramdim}$, which is the same structure as \eqref{eq:gradient_of_hamiltonian_of_discretized_system} in the proof of Proposition~\ref{thm:structure-preserving_neural_galerkin}.
	Thus, the conservation of $\Hsampled$ follows from analogous arguments as the ones of the proof of Proposition~\ref{thm:structure-preserving_neural_galerkin}.
\end{proof}

\begin{remark}
	The matrix function $\Jsampledblock{22}$ is obtained by exploiting the pointwise skew-symmetry of $\vecJ^{(22)}$ to ensure the pointwise skew-symmetry of $\Jsampledblock{22}$. However, this construction of $\Jsampledblock{22}$ cannot guarantee that \eqref{eq:structuredSampledSystem} with given $\param(t)$ and singular $\Msampled(\param(t))$ can be satisfied for a $\dot\param(t)$. The same issue applies to the special case of linear parametrizations.
\end{remark}

\section{Conserving quantities in time-discrete Neural Galerkin approximations}
\label{sec:modified_discretization}
We now consider the time discretization of the constrained or weighted Neural Galerkin systems, which is delicate because the nonlinear dependence of the parametrization $\ansatz(\param(t), \cdot)$ on the parameter $\param(t)$ means that quantities \eqref{eq:conserved_quantity_format} become nonlinear in $\param(t)$ and thus are not conserved by just applying Runge-Kutta integrators; see also the problem formulation in Section~\ref{sec:problem_statement}.
Building on literature of ODE integrators \cite[Chapter VII.2]{HairerW1996Solving}, we propose to use a nonlinear projection method that computes embeddings to conserve quantities in Neural Galerkin solutions in discrete time; see Figure~\ref{fig:manifold_illustration}b.
Importantly, the nonlinear projection approach is applicable with implicit and explicit time integration schemes. Explicit schemes are of especially great interest in the context of nonlinear parametrizations with Neural Galerkin schemes because they lead to linear least-squares regression problems at each time step whereas implicit schemes lead to non-convex optimization problems at each time step \cite{BrunaPV2022Neural,BermanP2023Randomized}.
We also discuss an alternative time discretization scheme which is implicit and based on discrete gradients.

\subsection{Neural Galerkin schemes with embeddings for conserving quantities in discrete time}
Consider a one-step time integrator applied to the constrained Neural Galerkin system \eqref{eq:modified_ng_ode} or the weighted Neural Galerkin system \eqref{eq:NonautonomousSemiDiscretizedSystem} with time-step size $\delta t > 0$. Such an integrator gives rise to a map $\vecPhi_{\delta t}: \Theta \to \Theta$ that takes a parameter vector $\param_k \in \Theta$ at time step $k$ and maps it onto $\param_{k + 1} \in \Theta$ at time step $k + 1$.
Let $\param_0$ be the parameter of the initial condition $\ansatz(\param_0, \cdot)$ and let $\param_1, \dots, \param_K \in \Theta$ be the parameters obtained with the time integrator at time steps $k = 1, \dots, K$, which lead to the time-discrete Neural Galerkin solution trajectory $\ansatz(\param_1, \cdot), \ansatz(\param_2, \cdot), \dots, \ansatz(\param_K, \cdot)$.
At each time step $k$, we want to ensure that the function $\ansatz(\param_{k + 1}, \cdot)$ at the next time step $k + 1$ is on the constrained manifold $\manifoldC$ if the approximate solution field function $\ansatz(\param_k, \cdot)$ at the current time step $k$ is on $\manifoldC$.
To achieve this, at each time step $k = 0, \dots, K-1$, we evaluate $\vecPhi_{\delta t}$ at $\param_k$  to compute $\tilde\param_{k + 1}$, which can correspond to a function $\ansatz(\tilde\param_{k + 1}, \cdot)$ that is outside of $\manifoldC$. We then follow \cite[Chapter IV]{HairerWL2006Geometric} and embed the function $\ansatz(\tilde\param_{k + 1}, \cdot)$ corresponding to the parameter $\tilde\param_{k + 1}$ as a function with parameter $\param_{k + 1}$ onto $\manifoldC$ by solving for $\param_{k + 1}$ via the nonlinear least-squares problem
\begin{align}
	\label{eq:manifold_minimization_problem}
	\min\limits_{\ansatz(\veceta, \cdot) \in \manifoldC} \frac{1}{2}\| \veceta - \tilde\param_{k + 1} \|^2_2\,.
\end{align}
The embedding step constrains the manifold on which we seek approximations to the PDE solution to $\manifoldC$, where quantities are conserved. The analogous procedure can be derived for embeddings onto $\hat\manifoldC$ defined in \eqref{eq:NonlinProj:SampledME},
\begin{equation}
	\label{eq:sampled_manifold_minimization_problem}
	\min\limits_{\ansatz(\veceta, \cdot) \in \hat\manifoldC} \frac{1}{2}\| \veceta - \tilde\param_{k + 1} \|^2_2\,.
\end{equation}

\subsection{Implicit time discretizations with discrete gradients}
Discretizing the second block equation in \eqref{eq:modified_ng_ode} via discrete gradients also leads to Neural Galerkin approximations that conserve quantities~\cite{Gonzalez1996Time,McLachlanQR1999Geometric}. The time integration via discrete gradients is implicit in time and thus can be computationally expensive in the context of Neural Galerkin schemes with nonlinear parametrizations because a potentially non-convex optimization problem has to be solved at each step; see \cite{BrunaPV2022Neural,BermanP2023Randomized}.
A discrete gradient \cite{Gonzalez1996Time,McLachlanQR1999Geometric} of a continuously differentiable function $H\colon \R^n\to\R$ is given by a continuous mapping $\overline{\nabla}H\colon \R^n\times\R^n\to\R^n$ which satisfies
\begin{equation}
	\label{eq:defining_properties_of_discrete_gradients}
	\overline{\nabla}H(\veceta,\veceta) = \nabla H(\veceta),\quad \overline{\nabla}H(\veceta,\veceta')^\T(\veceta'-\veceta) = H(\veceta')-H(\veceta)
\end{equation}
for all $\veceta,\veceta' \in\R^n$.
Based on discrete gradients, we discretize the second block equation $g(\param(t))^\T\dot{\param}(t)=0$ in \eqref{eq:modified_ng_ode} by replacing $\dot{\param}(t)$ by a finite difference approximation with time-step size $\delta t > 0$ and the gradients in \eqref{eq:CNG:GVecCont} by corresponding discrete gradients. 
For $k = 0, \dots, K - 1$, the resulting time-discrete equation is
\begin{equation}
	\label{eq:discreteGradientDiscretization}
	\overline{g}(\param_k,\param_{k+1})^\T \frac{\param_{k+1}-\param_k}{\Delta t} = 0,
\end{equation}
with $\overline{g}(\param_k,\param_{k+1}):=\begin{bmatrix}
		\overline{\nabla} \cq_1(\param_k,\param_{k+1}), \dots, \overline{\nabla} \cq_{\ncq}(\param_k,\param_{k+1})
	\end{bmatrix}$ using the discrete gradients $\overline{\nabla} \cq_i(\param_k,\param_{k+1})$ for $i = 1, \dots, \ncq$.
By using the second relation in \eqref{eq:defining_properties_of_discrete_gradients}, we infer that the $i$th row of \eqref{eq:discreteGradientDiscretization} is
\begin{equation*}
	0 = \overline{\nabla} \cq_i(\param_k,\param_{k+1})^\T \frac{\param_{k+1}-\param_k}{\Delta t} = \frac{\cq_i(\param_{k+1})-\cq_i(\param_k)}{\Delta t}
\end{equation*}
for $i=1,\ldots,\ncq$.
Consequently, \eqref{eq:discreteGradientDiscretization} is equivalent to
\begin{equation}
	\label{eq:time_discrete_preservation_of_conserved_quantities}
	\cq_1(\param_{k+1}) = \cq_1(\param_k),\quad \ldots,\quad \cq_{\ncq}(\param_{k+1}) = \cq_{\ncq}(\param_k),
\end{equation}
i.e., in this setting we do not have to construct the discrete gradients explicitly, but we may directly use \eqref{eq:time_discrete_preservation_of_conserved_quantities}.
The first equation in \eqref{eq:modified_ng_ode} can be discretized by any time discretization scheme.
Similarly, discrete gradients may be also used to ensure conservation in the context of the weighted Neural Galerkin systems \eqref{eq:NonautonomousSemiDiscretizedSystem}; see \cite[App.~C]{Schulze2023Energy-based_thesis} for more details.

\section{Computational aspects of Neural Galerkin with embeddings}
\label{sec:computational_aspects}
We now describe computational aspects of constrained Neural Galerkin schemes with nonlinear embeddings for conserving quantities.
In particular, we use the specific iterations as in \cite{HairerWL2006Geometric} to efficiently perform the embedding onto the constrained manifold. While we focus on constrained systems with embeddings, similar observation holds for the weighted Neural Galerkin schemes and time discretizations with discrete gradients.

\subsection{Constrained Neural Galerkin schemes represented as least-squares problems}
\label{subsec:least_squares}
While the sampled Neural Galerkin system \eqref{eq:Prelim:SampledNGODE} introduced in Section~\ref{sec:NG} can be numerically solved directly, the  corresponding linear least-squares problem is typically better conditioned. Recall that $\vecx_1, \dots, \vecx_{\nsamples}$ are the sample points used to obtain the sampled Neural Galerkin system \eqref{eq:Prelim:SampledNGODE}. If the solution field is scalar valued so that $m = 1$ and we use explicit Euler to discretize time with time-step size $\delta t > 0$, then the corresponding sampled least-squares problem at time steps $k = 0, \dots, K-1$ is
\begin{equation}\label{eq:Comp:Lsq}
	\min_{\delta\param_k \in \Theta} \|\hat{\vecA}(\param_k)\delta\param_k - \hat{\vecb}(\param_k)\|_2^2\,,
\end{equation}
where $\param_{k + 1} = \param_k + \delta t \delta\param_k$ and $\hat\vecA(\param_k) \in \mathbb{R}^{\nsamples \times \paramdim} $ and $\hat\vecb(\param_k) \in \R^{\nsamples}$ are
\begin{equation}\label{eq:Comp:Ab}
	\hat\vecA(\param_k)  = \begin{bmatrix}
		\nabla_{\param}\ansatz(\param_k, \vecx_1)^\T \\
		\vdots                                       \\
		\nabla_{\param}\ansatz(\param_k, \vecx_{\nsamples})^\T\end{bmatrix}\,,\quad
	\hat\vecb(\param_k)  = \begin{bmatrix}
		\vecf(\cdot, \ansatz(\param_k, \cdot))(\vecx_1)         \\
		\vdots                                                  \\
		\vecf(\cdot, \ansatz(\param_k, \cdot))(\vecx_\nsamples) \\
	\end{bmatrix}\,.
\end{equation}
The conservation constraints are added to \eqref{eq:Comp:Lsq} analogous to~\eqref{eq:modified_ng_ode}, which results in a least-squares problem with equality constraints
\begin{equation}
	\label{eq:Comp:LsqWithConstraints}
	\min_{\delta\param_k \in \Theta} \left\|\hat\vecA(\param_k)\delta\param_k - \hat\vecb(\param_k)\right\|_2^2  \, \quad \text{s.t. } \hat \vecg(\param_k)^\T \delta\param_k = 0,
\end{equation}
where $\hat \vecg$ is defined in \eqref{eq:TimeCont:SampledGVec}.
An analogous least-squares problem with equality constraints can be derived for solution fields with multiple outputs $m > 1$ and other time-integration schemes.

\subsection{Computing embeddings onto the constrained manifold}
To numerically solve~\eqref{eq:manifold_minimization_problem} with respect to the sampled manifold $\hat\manifoldC$, we introduce the Lagrangian function
\begin{equation}
	\label{eq:lagrangian}
	(\veceta, \veclambda) \mapsto \frac{1}{2}\|\veceta - \tilde{\param}_{k+1} \|_2^2 + \veclambda \cdot \manifoldFun(\veceta),
\end{equation}
where $\cdot$ denotes the Euclidean inner product and the differentiable function $\manifoldFun : \R^{\paramdim} \rightarrow \R^{\ncq}$ is defined as
\begin{align}
	\label{eq:manifoldFun}
	\manifoldFun(\veceta) =
	[
	\hat{\cq}_1(\ansatz(\veceta, \cdot)) - \hat{\cq}_1(\ansatz(\param(0), \cdot)), \dots,
	\hat{\cq}_\ncq(\ansatz(\veceta, \cdot)) - \hat{\cq}_\ncq(\ansatz(\param(0), \cdot))
	]^\T.
\end{align}
Thus,  $\manifoldFun(\veceta) = 0$ implies  that $\ansatz(\veceta, \cdot) \in \hat\manifoldC$.
The first-order optimality condition leads to the system of equations
\begin{subequations}
	\begin{align}
		\label{eq:projection_optimality_a}
		\veceta & = \parampreprojection + \manifoldFun'(\veceta)^\T \veclambda, \\
		\label{eq:projection_optimality_b}
		0       & = \manifoldFun(\veceta).
	\end{align}
\end{subequations}
We follow \cite[Section~IV.4]{HairerWL2006Geometric} and solve \eqref{eq:projection_optimality_a}--\eqref{eq:projection_optimality_b} via the iterations
\begin{subequations}
	\label{eq:newton_iteration}
	\begin{align}
		\Delta \veclambda^{(i)} & = -(\manifoldFun'(\parampreprojection)\manifoldFun'(\parampreprojection)^\T)^{-1}\manifoldFun(\parampreprojection + \manifoldFun'(\parampreprojection)^\T \veclambda^{(i)}), \\
		\veclambda^{(i+1)}      & = \veclambda^{(i)} + \Delta \veclambda^{(i)}\,,
	\end{align}
\end{subequations}
for $i = 1, 2, 3, \dots$.
Notice that the iterations are over the quantity $\veclambda$, which is of dimension $\ncq$. The number of quantities $\ncq$ is often much smaller than the number of parameters $\paramdim$.

\subsection{Neural Galerkin schemes with embeddings}

\begin{algorithm}[t]
	\caption{Neural Galerkin scheme with embeddings}
	\label{alg:embedded_ng}
	\begin{algorithmic}[1]
		\Procedure{NGEmbedding}{$\vecf$, $\delta t$, $K$, $\ansatz$, $\nsamples$, $\param_0$, $\cq_1, \dots, \cq_{\ncq}$, $\nsamplesmani$}
		\For{$k = 0, \dots, K - 1$}
		\State{Estimate $\hat\vecA(\param_k)$ and $\hat \vecb(\param_k)$ defined in \eqref{eq:Comp:Ab} using $\nsamples$ sample points.}
		\State{Estimate $\hat\vecg(\param_k)$ using $\nsamplesmani$ sample points and quantities $\cq_1, \dots, \cq_{\ncq}$.}
		\State{Compute $\delta \param_k$ as solution of~\eqref{eq:Comp:LsqWithConstraints}.}
		\State{Set $\tilde \param_{k+1} = \param_{k} + \delta t \delta \param_k$.}
		\State{Iterate~\eqref{eq:newton_iteration} to compute $\param_{k+1}$.}
		\EndFor
		\State Return trajectory $\ansatz(\param_0, \cdot), \ansatz(\param_1, \cdot), \dots, \ansatz(\param_K, \cdot)$
		\EndProcedure
	\end{algorithmic}
\end{algorithm}
Algorithm~\ref{alg:embedded_ng} summarizes the Neural Galerkin schemes with embeddings when explicit Euler is used to discretize time. We stress that other time integration schemes can be used in an analogous way.

\subsubsection{Description of algorithm}
The inputs to Algorithm~\ref{alg:embedded_ng} are the right-hand side function $\vecf$ of the PDE \eqref{eq:general_pde}, time-step size $\delta t$, number of time steps $K$, parametrization $\ansatz$, number of sampling points $\nsamples$, parameter $\param_0$ corresponding to the initial condition $\ansatz(\param_0, \cdot)$, the quantities $\cq_1, \dots, \cq_{\ncq}$ and the number of sampling points $\nsamplesmani$ for the sampled quantities $\hat\cq_1, \dots, \hat\cq_{\ncq}$.
The algorithm iterates over the time steps $k = 0, \dots, K - 1$. In each iteration, the least-squares problem \eqref{eq:Comp:Ab} and the sampled gradients \eqref{eq:TimeCont:SampledGVec} are assembled.
Notice that the sampled 
$\hat{\vecg}$ is used, which is based on the sampled quantities $\hat\cq_1, \dots, \hat\cq_{\ncq}$ with $\nsamplesmani$ sampling points.
Then, the update $\delta \param_k$ is computed by solving the least-squares problem \eqref{eq:Comp:LsqWithConstraints}.
In line 7, the vector $\tilde{\param}_{k + 1} = \param_k + \delta t \delta \param_k$ is projected onto $\param_{k + 1}$ so that the solution field $\ansatz(\param_{k + 1}, \cdot)$ lies on the constrained manifold.
After $K$ time steps, the trajectory of Neural Galerkin approximations $\ansatz(\param_0, \cdot), \dots, \ansatz(\param_K)$ is returned.

\subsubsection{Computational costs}
The costs of a time step of Algorithm~\ref{alg:embedded_ng} is dominated by estimating the sampled  $\hat\vecA(\param_k)$, $\hat \vecb(\param_k)$, and $\hat{\vecg}(\param_k)$ and subsequently solving the least-squares problem \eqref{eq:Comp:LsqWithConstraints} to obtain $\delta \param_k$.
The costs of computing the sampled $\hat\vecA(\param_k), \hat \vecb(\param_k)$ scale with the number of sampling points $\nsamples$. For each sample, the parametrization $\ansatz$, its gradient $\nabla_{\param}\ansatz$, and the right-hand side function $\vecf$ are evaluated. The costs of these evaluations depend on the parametrization. If a fully connected deep network with $\ell$ layers and $p$ nodes is used, then the costs scale as $\mathcal{O}(\ell p^2)$.
The costs of computing $\hat{\vecg}(\param_k)$ are dominated by the quadrature in \eqref{eq:conserved_quantity_monte_carlo}, which depends on the number of sampling points $\nsamplesmani$. The functions $\kappa_1, \dots, \kappa_{\ncq}$ that define $\cq_1, \dots, \cq_{\ncq}$ as shown in \eqref{eq:conserved_quantity_monte_carlo} are typically cheap to evaluate.
The costs of solving the least-squares problem \eqref{eq:conserved_quantity_format}  scales as $\mathcal{O}(\nsamples \paramdim^2)$ if the generalized QR decomposition~\cite{AndersonBD1992Generalized} is used, which is implemented in the LAPACK routine \texttt{DGGLSE}. We refer to \cite{BermanP2023Randomized} for work on reducing the costs per time step.
The projection \eqref{eq:newton_iteration} computed in line 7 of Algorithm~\ref{alg:embedded_ng} typically requires only few iterations and thus incurs negligible costs in our numerical experiments.

\section{Numerical experiments}
\label{sec:numerical_experiments}
We now demonstrate Neural Galerkin with embeddings with numerical experiments. First, we consider the inviscid Burgers' equation~\eqref{eq:burgers_equation_hamiltonian_structure} to demonstrate the interplay between the constrained Neural Galerkin system and the nonlinear embedding.
Second, we illustrate with the acoustic wave equation~\eqref{eq:wave_equation_hamiltonian_structure} how the number of sampling points $\nsamplesmani$ influences the preservation of the Hamiltonian \eqref{eq:wave_equation_hamiltonian_structure}.
Third, we consider the shallow water equations in a two-dimensional spatial domain with the total energy as conserved quantity.
The implementation is based on \texttt{jax}, which is a Python library for automatic differentiation~\cite{jax2018github}. All floating point computations are carried out in double precision. The implementation is available online  \url{github.com/Algopaul/ng_embeddings}.

\subsection{Experimental setup}

In each experiment, we compare approximations to reference solutions $\vecval{u}^{\text{ref}}$ obtained with spectral methods. Time is discretized in all examples with the explicit fourth-order Runge-Kutta method~\cite[Table 1.2, l]{HairerWN1993Solving}. For a solution trajectory with parameters $\{\param_k\}_{k=0}^{K}$, $K \in \N$, at time step $\{t_k\}_{k=0}^{K}$, we compute the relative error of the solution field as
\begin{equation}
	\label{eq:rel_err_def}
	E_r(t_k) = \sum\limits_{i=1}^{n_E}\|\ansatz(\param_k, \vecx^{\text{test}}_i) - \vecval{u}^{\text{ref}}(t_k, \vecx^{\text{test}}_i)\| \Big/ \sum\limits_{i=1}^{n_E}\|\vecval{u}^{\text{ref}}(t_k, \vecx^{\text{test}}_i)\|,
\end{equation}
using $n_E$ equidistantly sampled test points $\vecx^{\text{test}}_1, \dots, \vecx^{\text{test}}_{n_E}$ in the respective domains.
We also report the error in conserving the quantities of interest $t \mapsto \|\hat{\cq}^{\text{test}}(t)-\hat{\cq}^{\text{test}}(0)\|$ for the different solution trajectories, where we use $n_E$ equidistantly sampled test points in the respective domain to estimate the quantity, which we denote as $\hat{\cq}^{\text{test}}$.
We ensure that the $n_E$ points used during the evaluation are different from the $\nsamples$ points used to construct the least-squares problem~\eqref{eq:Comp:Ab} and the $\nsamplesmani$ points used to estimate the quantities during time integration.

\subsection{Burger's equation with conservation of mass} \label{sec:burgers_equation}
\begin{figure}[t]
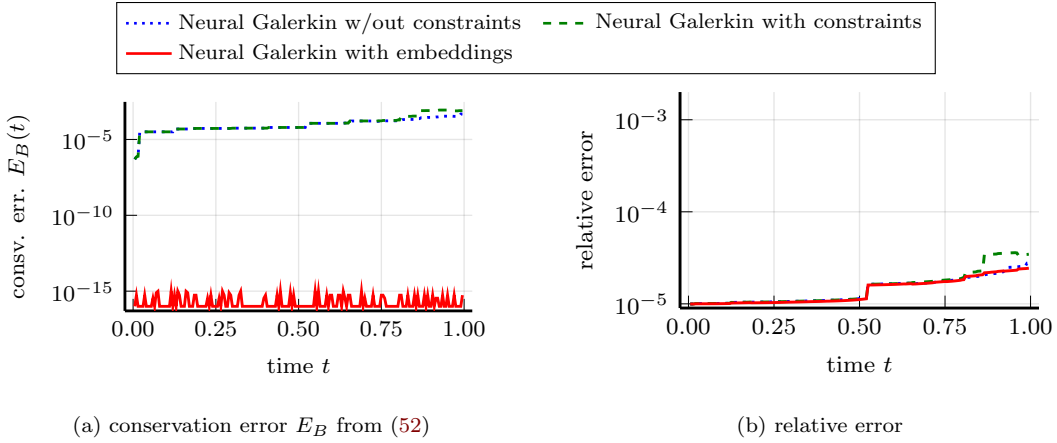

	\centering
	\ref{burgerlegend} \\
	\begin{tabular}{cc}
		\resizebox{0.48\columnwidth}{!}{\input{./PlotSources/burgers_conservation_error_mlp.tikz}}
		                                                                                                         & \resizebox{0.48\columnwidth}{!}{\input{./PlotSources/burgers_output_error_mlp.tikz}}
		\\
		\scriptsize{(a) conservation error $\burgersconservationerror$ from~\eqref{eq:burgersconservationerror}} &
		\scriptsize{(b) relative error}                                                                                                                                                                 \\
	\end{tabular}
	\caption{Burgers' equation: Nonlinear parametrization such as deep networks used in Neural Galerkin schemes imply that even linear quantities become nonlinear in the parameters. Thus, adding a constraint to the Neural Galerkin scheme is insufficient for conserving quantities with explicit Runge-Kutta schemes. In contrast, the proposed Neural Galerkin scheme with embeddings computes projections onto constrained manifolds at each time step to conserve quantities even in discrete time.}
	\label{fig:comparison_linear_nonlinear}
\end{figure}

We consider the inviscid Burgers' equation \eqref{eq:burgers_equation_hamiltonian_structure} in the spatial domain $[-1, 1)$ with periodic boundary conditions and conserved quantity $\cq_{\text{mass}} = \int_{-1}^1 u \, \dd x$.
We parametrize the solution field with a fully connected feed-forward deep network that has three hidden layers of width ten with sinusoidal activation functions. The first layer imposes periodicity as in~\cite{BermanP2023Randomized}.
The output layer is linear. It is important to note that the conserved quantity is linear in the solution function $u$ but nonlinear in the parameter $\param(t)$ of the parametrization.
We use $n_S=200$ equidistant sample points to construct the least-squares problem~\eqref{eq:Comp:Lsq} and $n_M=200$ equidistant points to estimate the conserved quantity and its derivative~\eqref{eq:TimeCont:SampledGVec}.
The time-step size is $\delta t = 5 \cdot 10^{-3}$.
We now compare Neural Galerkin without constraints \eqref{eq:Prelim:SampledNGODE}, Neural Galerkin with constraints \eqref{eq:modified_ng_ode_sampled} that conserves quantities in continuous time, and Neural Galerkin with embeddings that is described in Algorithm~\ref{alg:embedded_ng} and that combines \eqref{eq:modified_ng_ode_sampled} with the embeddings \eqref{eq:sampled_manifold_minimization_problem} to conserve quantities also in discrete time.  Figure~\ref{fig:comparison_linear_nonlinear}a shows the error in conserving the quantity $\hat{\cq}^{\text{test}}_{\text{mass}}$, estimated with $n_E = 400$, which we define as
\begin{align}
	\label{eq:burgersconservationerror}
	\burgersconservationerror(t) = |\hat{\cq}^{\text{test}}_{\text{mass}}(t) -\hat{\cq}^{\text{test}}_{\text{mass}}(0)|.
\end{align}
Because we use a deep network which is a nonlinear parametrization, the quantity is not conserved when just adding a constraint. In contrast, the proposed Neural Galerkin scheme with embeddings conserves the sampled quantity because it performs an explicit embedding after each time step. Figure~\ref{fig:comparison_linear_nonlinear}b shows the relative error computed using~\eqref{eq:rel_err_def} with $n_E=400$. For the reference solution, we have used 200 Fourier modes and a time step size of $10^{-3}$.

\subsection{Acoustic wave equation with conservation of Hamiltonian} \label{sec:wave_equation}
\begin{figure}[t]
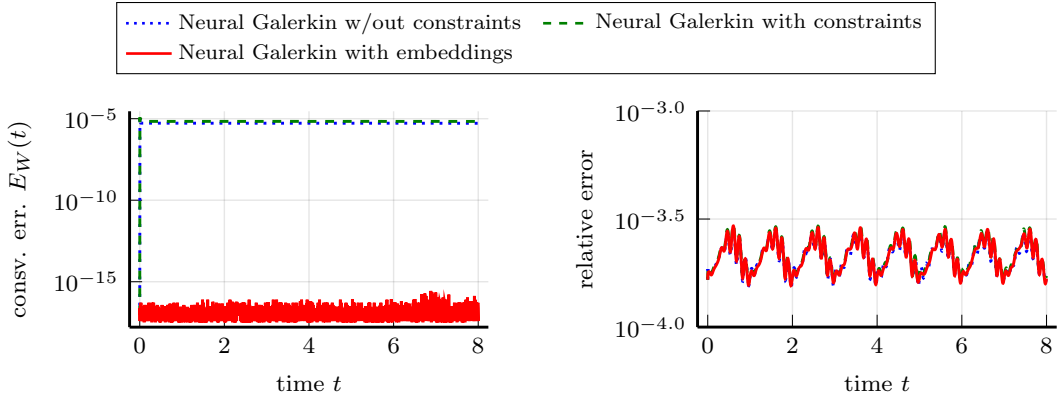

	\begin{center}
		\ref{hamilwaveoned} \\
		\begin{tabular}{cc}
			\resizebox{0.48\columnwidth}{!}{\input{./PlotSources/hamil_wave_1d_hamiltonian_error.tikz}} &
			\resizebox{0.48\columnwidth}{!}{\input{./PlotSources/hamil_wave_1d_output_error.tikz}}        \\
			\scriptsize{(a) error in conserving Hamiltonian}                                           &
			\scriptsize{(b) relative error in solution fields}                                           \\
		\end{tabular}
	\end{center}
	\caption{Acoustic wave equation: The proposed Neural Galerkin scheme with embeddings preserves the Hamiltonian at test points to machine precision.}
	\label{fig:hamil_wave_1d_errors}
\end{figure}

\begin{figure}[t]
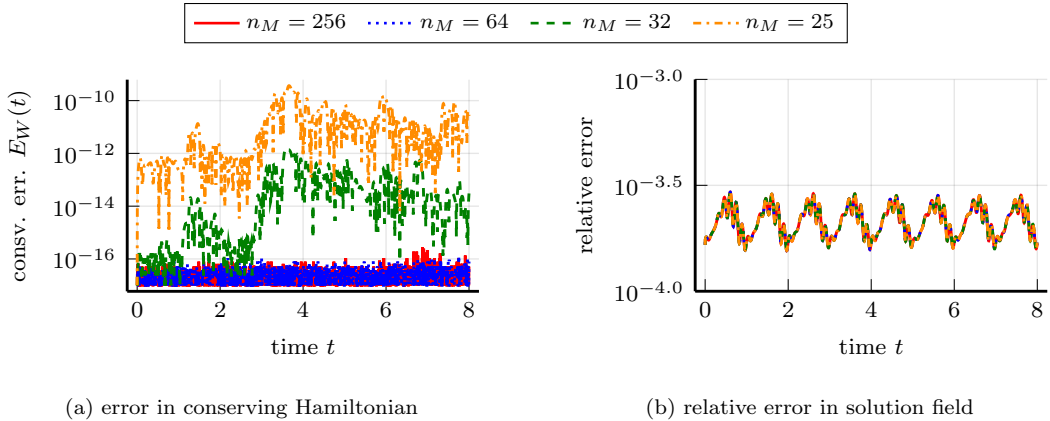

	\begin{center}
		\ref{hamilwaveonedsubsamples}
		\begin{tabular}{cc}
			\resizebox{0.48\columnwidth}{!}{
				\input{./PlotSources/hamil_wave_1d_rel_errs_time_subsample.tikz
			}} &
			\resizebox{0.48\columnwidth}{!}{
				\input{./PlotSources/hamil_wave_1d_ham_errs_time_subsample.tikz
			}} \\
			\scriptsize{(a) error in conserving Hamiltonian} &
			\scriptsize{(b) relative error in solution field}
		\end{tabular}
		\caption{Acoustic wave equation: The proposed Neural Galerkin scheme with embeddings preserves the Hamiltonian at a fine grid of test points even when a low number of sampling points $\nsamplesmani$ is used.}
		\label{fig:wave1d_subsampling}
	\end{center}
\end{figure}

We present results for the acoustic wave equation \eqref{eq:wave_equation_hamiltonian_structure} with its Hamiltonian $H_{\text{wave}}$ as conserved quantity.
We have periodic boundary conditions on the spatial domain $\spacedomain = [-1, 1)$ with end time $T = 8$ and initial condition $\rho(0, x) = \mathrm{e}^{-9x^2}, v(0, x) = 0$. We set $\rhoref=c=1$.
The parametrization is a fully connected feed-forward network with two hidden layers of width ten followed by a hidden layer of width 20, each with sinusoidal activation functions, a periodic input layer as in~\cite{BermanP2023Randomized} of width ten, and a linear output layer of width two to account for $\rho$ and $v$.
We use $n_S=256$ sampling points distributed equidistantly in $\spacedomain$ to assemble the least-squares problem~\eqref{eq:Comp:Lsq} and the same $n_M=256$ sampling points to estimate the Hamiltonian and its derivative~\eqref{eq:TimeCont:SampledGVec}.
We discretize time using fourth-order Runge-Kutta with a time-step size of $2^{-8} \approx 4 \times 10^{-3}$.
For the reference solution, we have used 256 Fourier modes and fourth-order explicit Runge-Kutta with a time-step of $10^{-3}$.
Figure~\ref{fig:hamil_wave_1d_errors}a shows the error in conserving the Hamiltonian. We denote the error in $\hat{H}_{\text{wave}}$ by
\begin{equation}
	\label{eq:ham_cons_quant_error}
	E_{W}(t) = |\hat{H}^{\text{test}}_{\text{wave}}(t) - \hat{H}^{\text{test}}_{\text{wave}}(0)|,
\end{equation}
where we set $n_E = 512$. The results show that our scheme based on embeddings preserves the sampled Hamiltonian at test points whereas the other schemes lead to solutions with large variations in the Hamiltonian. Figure~\ref{fig:hamil_wave_1d_errors}b shows the relative error~\eqref{eq:rel_err_def} estimated with $n_E=512$ samples. The relative error of all three schemes is comparable in this example.
In Figure~\ref{fig:wave1d_subsampling}, we demonstrate that we can use fewer samples $\nsamplesmani$ for the estimation of $\hat{H}_{\text{wave}}$ during the time integration and still achieve conservation at test points $\hat{H}_{\text{wave}}^{\text{test}}$ to machine precision. In fact, in Figure~\ref{fig:wave1d_subsampling}, it is shown the sampled Hamiltonian at the test points is conserved to machine precision even when choosing $\nsamplesmani=64$ and still to $10^{-10}$ for $\nsamplesmani=25$.

\subsection{Shallow water waves in two spatial dimensions with energy conservation}\label{sec:shallow_water}
We now consider the shallow water equations in a two-dimensional domain with periodic boundary conditions. We follow a similar setup as the one introduced in~\cite[Section 8]{HesthavenPR2022Rank-adaptive}.

\subsubsection{Setup}
The governing equations are
\begin{align}
	\label{eq:shallow_water2d}
	\begin{split}
		\partial_t \sweheight  + \nabla \cdot (\sweheight \nabla \swephi) &= 0\quad                 \text{ in } \swetimedomain\times\swespacedomain  , \\
		\partial_t \swephi + \frac{1}{2} {\left\| \nabla\swephi \right\|}^2 + \sweheight &= 0 \quad  \text{ in } \swetimedomain\times\swespacedomain, \\
	\end{split}
\end{align}
with initial conditions
\begin{align}
	\label{eq:shallow_wated2d_init}
	\begin{split}
		\sweheight(0, \vecx) &= 1+0.33 \mathrm{e}^{-1.7 {\left\|\vecx\right\|}^2} \quad             \text{ in } \swespacedomain,                       \\
		\swephi(0, \vecx) &= 0 \quad                                                                 \text{ in } \swespacedomain,                       \\
	\end{split}
\end{align}
where $\swespacedomain = [-4, 4)^2$, and $\sweheight, \swephi: \swetimedomain\times\swespacedomain   \rightarrow \R$ denote the height and the potential field of the fluid.
A conserved quantity of~\eqref{eq:shallow_water2d} is the energy
\begin{align}
	\label{eq:swe_conserved_quantity}
	\swehamiltonian(\sweheight, \swephi) = \frac{1}{2}
	\int_{\swespacedomain} \sweheight {\left\| \nabla\swephi \right\|}^2 + h^2 \dd \vecx.
\end{align}
We scale the equations above such that $\sweheight$ and $\swephi$ are close in magnitude and centered around $0$ to avoid numerical issues,
\begin{align}
	\label{eq:shallow_water2d_modified}
	\begin{split}
		\partial_t \widetilde \sweheight  + \nabla \cdot ((\widetilde \sweheight +1)\nabla \swephi) = 0,                       & \quad \text{ in } \swetimedomain\times\swespacedomain, \\
		\partial_t \widetilde \swephi + \frac{1}{2} {\left\| \nabla\widetilde \swephi \right\|}^2 + \widetilde \sweheight = 0, & \quad \text{ in } \swetimedomain\times\swespacedomain\,. \\
	\end{split}
\end{align}
The corresponding initial conditions are
\begin{align}
	\label{eq:shallow_water2d_modified_init}
	\begin{split}
		\widetilde \sweheight(0, \vecx) = 0.33 \mathrm{e}^{-1.7 {\left\|\vecx\right\|}^2},                                     & \quad \text{ in } \swespacedomain,                       \\
		\widetilde \swephi(0, \vecx) = 0,                                                                                      & \quad\text{ in } \swespacedomain,                       \\
	\end{split}
\end{align}
where $\widetilde \sweheight = \sweheight - 1$ and $\widetilde \swephi = \swephi + t$. Note that in the second equation, we did not add the constant $1$ to prevent $\widetilde \swephi$ from diverging numerically.
Our reference solution is obtained with a spectral method with $300$ degrees of freedom in each spatial dimension and a time-step size $10^{-3}$.

\subsubsection{Results}
\begin{figure}[t]
	\begin{center}
		\begin{tabular}{ccc}
			\includegraphics[height=0.25\textwidth]{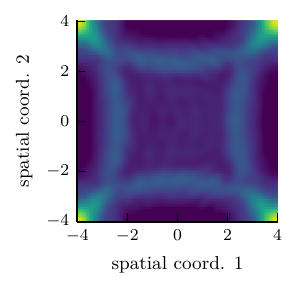}     &
			\includegraphics[height=0.25\textwidth]{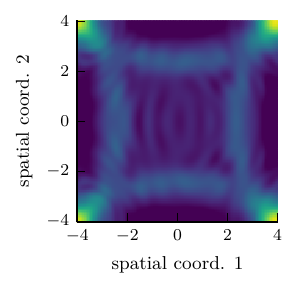}            &
			\includegraphics[height=0.25\textwidth]{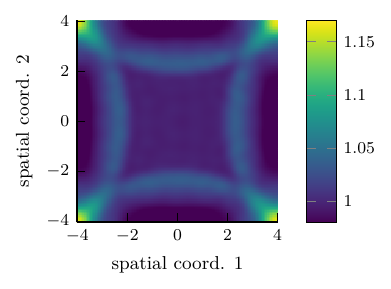}           \\
			\includegraphics[height=0.25\textwidth]{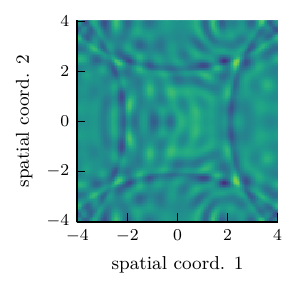} &
			\includegraphics[height=0.25\textwidth]{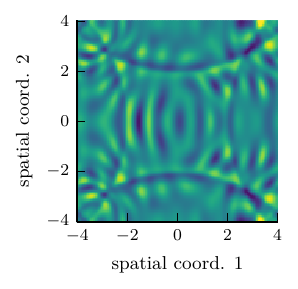}        &
			\includegraphics[height=0.26\textwidth]{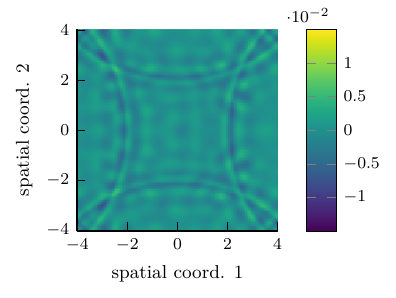}       \\
			\scriptsize{(a) Neural Galerkin}                                                           &
			\scriptsize{(b) Neural Galerkin with constraints}                                          &
			\scriptsize{(c) Neural Galerkin with embeddings}
		\end{tabular}
	\end{center}
	\caption{Shallow water: The proposed Neural Galerkin scheme with embeddings conserves energy over time and avoids oscillations in the solution field, which represents the height of the fluid. Time is $t = 5$.}
	\label{fig:swe2d_5sec_heatmap}
\end{figure}

\begin{figure}[t]
	\begin{center}
		\begin{tabular}{ccc}
			\includegraphics[height=0.25\textwidth]{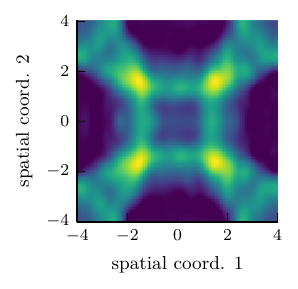}     &
			\includegraphics[height=0.25\textwidth]{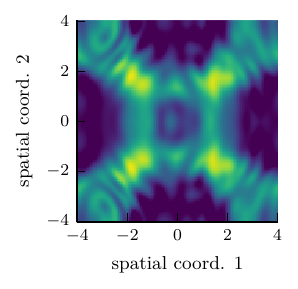}            &
			\includegraphics[height=0.25\textwidth]{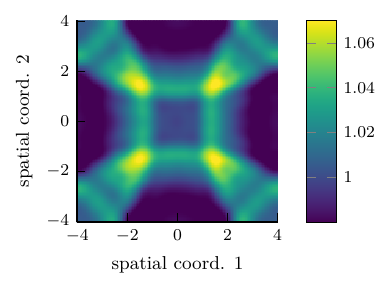}           \\
			\includegraphics[height=0.25\textwidth]{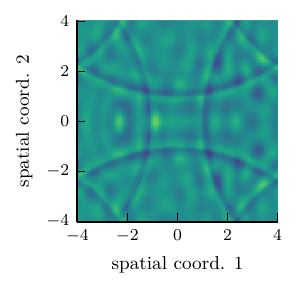} &
			\includegraphics[height=0.25\textwidth]{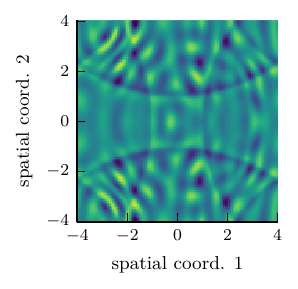}        &
			\includegraphics[height=0.26\textwidth]{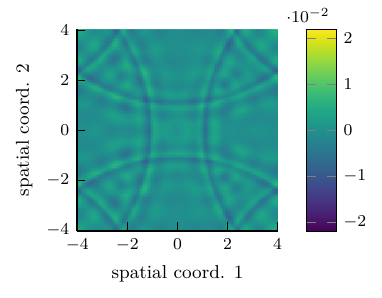}       \\
			\scriptsize{(a) Neural Galerkin}                                                            &
			\scriptsize{(b) Neural Galerkin with constraints}                                           &
			\scriptsize{(c) Neural Galerkin with embeddings}
		\end{tabular}
	\end{center}
	\caption{Shallow water: The plots of the solution field of $h$ obtained with Neural Galerkin and embeddings indicates that the conservation of energy avoids the oscillations that are present in the approximations obtained with the other schemes that ignore energy conservation. Time is $t = 6$.}
	\label{fig:swe2d_final_heatmap}
\end{figure}

\begin{figure}[t]
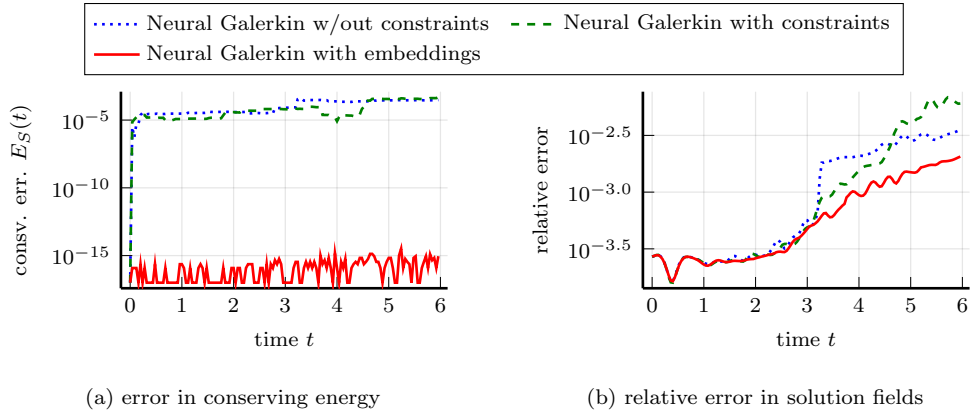

	\begin{center}
		\ref{swelegendvariants} \\
		\begin{tabular}{cc}
			\input{./PlotSources/swe2d_ham_errs_time.tikz} & \input{./PlotSources/swe2d_rel_errs_time.tikz} \\
			\scriptsize{(a) error in conserving energy}    &
			\scriptsize{(b) relative error in solution fields}                                              \\
		\end{tabular}
	\end{center}
	\caption{Shallow water: The proposed Neural Galerkin scheme with embeddings conserves the sampled energy at test points up to machine precision.}
	\label{fig:swe2d_error_hamiltonian}
\end{figure}

\begin{figure}[t]
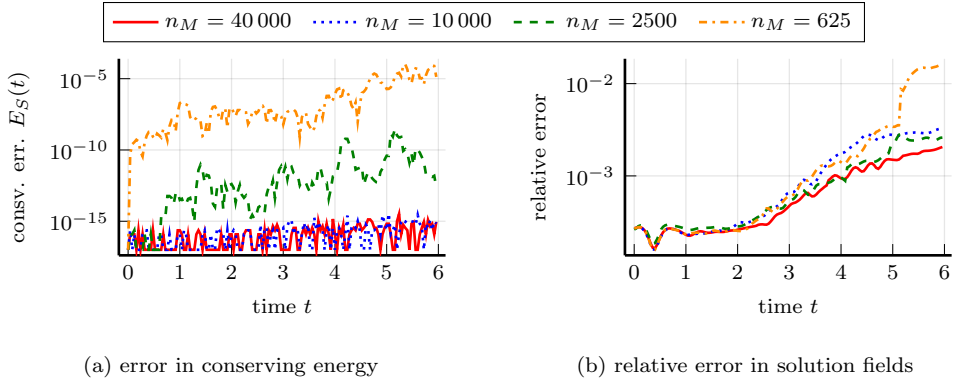

	\begin{center}
		\ref{swelegendsub} \\
		\begin{tabular}{cc}
			\input{./PlotSources/swe2d_subham_errs_time.tikz} & \input{./PlotSources/swe2d_subrel_errs_time.tikz} \\
			\scriptsize{(a) error in conserving energy}       &
			\scriptsize{(b) relative error in solution fields}                                                    \\
		\end{tabular}
	\end{center}
	\caption{Shallow water: The energy is conserved up to machine precision at the test points when $\nsamplesmani = 10\,000$ samples (100 in each spatial direction) are used and increases up to $10^{-5}$ when only $\nsamplesmani = 625$ samples are used (25 in each direction).}
	\label{fig:swe2d_error_hamiltonian_subsample}
\end{figure}

The parametrization is a deep network with the same structure as in Section~\ref{sec:wave_equation}. The time-step size is $\delta t = 2 \times 10^{-3}$ and the time discretization scheme is fourth-order explicit Runge-Kutta. We take 200 points in each spatial direction to assemble the least-squares problem~\eqref{eq:Comp:Lsq} and for estimating the conserved quantity~\eqref{eq:swe_conserved_quantity} and its derivative~\eqref{eq:TimeCont:SampledGVec}. We compare Neural Galerkin without constraints, Neural Galerkin with constraints, and the proposed scheme with embeddings.
Figure~\ref{fig:swe2d_5sec_heatmap} and~\ref{fig:swe2d_final_heatmap} show the solution field and the point-wise error of the fluid height $h$ at times $t=5$ and $t=6$, respectively. The approximation obtained with embeddings avoids the oscillations that are present in the approximations obtained with other schemes that ignore conservation of quantities.
The error in conserving the energy is shown in Figure~\ref{fig:swe2d_error_hamiltonian}a, which we measure as $E_S(t) = |\hat{I}^{\text{test}}_{(e)}(t)-\hat{I}^{\text{test}}_{(e)}(0)|$, for $n_E = 90\,000$ equidistantly sampled test points.
Figure~\ref{fig:swe2d_error_hamiltonian}b shows the relative error \eqref{eq:rel_err_def} with $n_E=90\,000$ test samples. The Neural Galerkin scheme with embeddings achieves the lowest relative error in the solution fields, which is in agreement with the plots shown in Figure~\ref{fig:swe2d_5sec_heatmap} and~\ref{fig:swe2d_final_heatmap} where the variants without energy conservation lead to oscillations in the approximate solution fields.
In Figure~\ref{fig:swe2d_error_hamiltonian_subsample} we study the effect of the the number of samples $\nsamplesmani$ used to estimate the energy.
The results show that with 100 points in each direction so that $\nsamplesmani = 10\,000$, the Neural Galerkin scheme with embeddings conserves the energy at the test points up to machine precision. The error in the energy conservation grows to about $10^{-5}$ only when there are only 25 sampling points in each spatial direction, which corresponds to a total of $\nsamplesmani = 625$ points.
Figure~\ref{fig:swe2d_error_hamiltonian_subsample}b shows how the relative error \eqref{eq:rel_err_def} in the solution fields depends on the number of sampling points $\nsamplesmani$. It can be seen that conserving the energy more accurately with more sampling points leads to lower relative errors in the solution fields in this example. If only $25$ sampling points in each spatial dimension are used, then the error increases by about one order of magnitude compared to when the energy is conserved to machine precision.

\section{Conclusions}
\label{sec:conclusion}
Preserving structure and conserving quantities in nonlinear approximations of PDE solutions is delicate because the nonlinear dependence on the parameter implies that even linear quantities in the solution fields can become nonlinear in the parameters.
Thus, just adding constraints to time-continuous formulations is insufficient with standard Runge-Kutta time integrators. While one can resort to implicit methods such as discrete gradients, it is desirable to have explicit time integration schemes when nonlinear parametrizations are used because they require solving systems of linear equations at each time step in typical cases rather than systems of nonlinear equations as with implicit integrators.
The proposed Neural Galerkin schemes compute explicit embeddings on manifolds of parametrizations that conserve quantities, which can be computed efficiently and are applicable with explicit time integrators and generic nonlinear parametrizations such as deep networks.

\printbibliography

\end{document}